\input amstex
\documentstyle{amsppt}
%\magnification=\magstep1
 \NoRunningHeads
 \loadbold
\topmatter

\NoBlackBoxes
\title
 Classification  of rank-one actions  via the cutting-and-stacking parameters
 \endtitle

 %\comment

\author
Alexandre  I. Danilenko  and Mykyta I. Vieprik
\endauthor

\address
Faculty of Mathematics and Computer Science,
Nicolaus Copernicus University,
Chopin street 12/18, 87-100 Toru{\'n}, Poland
\newline\indent
 B. Verkin Institute for Low Temperature Physics and Engineering
of the  National Academy of Sciences of Ukraine,
47 Nauky Ave.,
 Kharkiv, 61164, Ukraine
 \newline\indent
 Mathematical Institute of the Polish Academy of Sciences,
8  {\'S}niadeckich  street,
 Warsaw, 00-656,
 Poland
\endaddress
\email            alexandre.danilenko\@gmail.com
\endemail

\address
King’s College London, Department of Mathematics, Strand, London, WC2R 2LS, United Kingdom
 \newline\indent
 Mathematical Institute of the Polish Academy of Sciences,
 8 {\'S}niadeckich street,
 Warsaw, 00-656,
 Poland
\endaddress
\email
nikita.veprik\@gmail.com
\endemail

%\endcomment

\abstract
Let $G$ be a discrete countable  infinite group.
Let $T$ and $\widetilde T$ be two rank-one $\sigma$-finite measure preserving actions of $G$ and let $\Cal T$ and $\widetilde {\Cal T}$ be the cutting-and-stacking parameters that determine $T$ and $\widetilde T$ respectively.
We find necessary and sufficient conditions on $\Cal T$ and $\widetilde{\Cal T}$ under which   $T$ and $\widetilde T$ are
isomorphic.
% and  (ii) $\widetilde T$ is a factor of $T$.
We also show that the isomorphism equivalence relation is a $G_\delta$-subset
in the Cartesian square of the set of all admissible  parameters $\Cal T$ endowed with the natural Polish topology.
If $G$ is amenable and $T$ and $\widetilde T$ are finite measure   preserving then we also find 
necessary and sufficient conditioins on $\Cal T$ and $\widetilde {\Cal T}$ under which 
$\widetilde T$ is a  factor of $T$.
\endabstract

%\comment

\thanks
This work  was  supported in part by the 
``Long-term program of
support of the Ukrainian research teams at the 
Polish Academy of Sciences carried out in
collaboration with the U.S. National Academy of
Sciences with the financial support of external
partners''. 
 \endthanks

% \endcomment

 
\endtopmatter

\document

\head{0. Introduction}\endhead

Classification of ergodic dynamical systems up to isomorphism is a central problem of ergodic theory. 
However, despite some progress achieved  in the spectral classification of the ergodic transformations with discrete spectrum (Halmos-von Neumann theorem) or the  classification of the Bernoulli shifts via the Kolmogorov-Sinai entropy (Ornstein theorem), etc., it was  shown  rigorously  that no classification in a reasonable sense exists
for the entire class of ergodic systems.
We refer the reader to the survey \cite{Fo2} for more information and references to relevant ``non-classification'' works. 
This research direction was summed up with  a remarkable result by Foreman--Rudolph--Weiss: if the set of  
ergodic transformation  $\Cal E$ of a Lebesgue space $(X,\goth B,\mu)$ is endowed with the standard (Polish) weak operator topology
then the isomorphism equivalence relation {\bf Iso} on $\Cal E$, which is a subset of $\Cal E\times\Cal E$,
is not Borel  \cite{FoRuWe}. 
The restriction of {\bf Iso} to the subset of weakly mixing transformations
is not Borel either \cite{Ku}.
However, it was shown in \cite{FoRuWe, Theorem~51} that the restriction of {\bf Iso} to the subset $\Cal R_1\subset\Cal E$ of rank-one transformations of $(X,\goth B,\mu)$ is Borel.
It should be noted that $\Cal R_1$ is a dense $G_\delta$-subset of $\Cal E$.
Hence, $\Cal R_1$ is Polish in the induced topology.
In view of  \cite{FoRuWe, Theorem~51}, Foreman, Rudolph and Weiss stated a problem:
\roster
\item"$\circ$" {\it  to find a good explicit method
of checking when two rank-one transformations are isomorphic.}
\endroster
There exist many different ways to define rank-one transformations (see \cite{Fe} and references therein).
One of them is  the technique of cutting-and-stacking  with a single tower on every step 
of this inductive construction.
Then each rank-one transformation is completely determined by 
the underlying  cutting-and-stacking parameters: a sequence of cuts and a sequence of spacer mappings.
Therefore, we can
reformulate the aforementioned problem by Foreman--Rudolph--Weiss  as following:
\roster
\item"$\bullet$"
{\it to find necessary and sufficient conditions under which two families  of cut\-ting-and-stacking parameters determine isomorphic rank-one transformations.}
\endroster
It   is solved in the present work.
In fact, we solve this classification problem in a much more general setting of rank-one $\sigma$-finite measure preserving actions of arbitrary countable infinite discrete groups. 
The ``rank-one'' here means the rank one along a  sequence of finite subsets in the group.
If the invariant measure is finite then this sequence is necessarily F{\o}lner and the group is amenable.
If the group is $\Bbb Z$ then the rank-one finite measure preserving $\Bbb Z$-actions according to our definition are exactly
the  funny rank-one transformations (see \cite{Fe}).

 %for the rank-one probability preserving actions of arbitrary discrete countable amenable groups. 
%The ``rank-one'' means the rank one along a  F{\o}lner sequence in the group.
%In the case, where the acting group is $\Bbb Z$, this includes classification of 
 %funny rank-one transformations (see \cite{Fe}) along a fixed F{\o}lner sequence.

It is convenient to state and prove the main results in the language of $(C,F)$-systems.
This  is an algebraic version of the above mentioned geometric cutting-and-stacking construction.
It was introduced in \cite{dJ2} and \cite{Da1} in similar but non-equivalent ways.
We use below a  more general version of the $(C,F)$-construction from \cite{Da2} which embraces the
earlier versions from \cite{dJ2} and \cite{Da1} as particular cases. 
Each rank-one action is isomorphic to a $(C,F)$-action and the converse is also true (\cite{Da2}, \cite{DaVi}).
Each $(C,F)$-action of a group $G$ is determined by a certain sequence $(C_n, F_{n-1})_{n=1}^\infty$
of finite subsets $C_n$ and $F_{n-1}$ in $G$.
A pair $(C_n,F_{n-1})$  is simply an encoded information about how the copies of the $(n-1)$-th tower are located inside the $n$-th tower.

\proclaim{Theorem  A} Let $G$ be a countable infinite discrete group.
Let $T$ and $\widetilde T$ be two $\sigma$-finite measure preserving  $G$-actions associated with
$(C,F)$-sequences $(C_n,F_{n-1})_{n=1}^\infty$ and 
$(\widetilde C_n,\widetilde F_{n-1})_{n=1}^\infty$ respectively.
%Let $(\epsilon_n)_{n=1}^\infty$ be a summable sequence of positive reals.
Then $T$ and $\widetilde T$ are (measure theoretically) isomorphic if and only if there exist a sequence
$$
0=k_0=l_0=k_1<l_1<k_2<l_2<\cdots
$$
 of non-negative integers and subsets $J_n\subset F_{k_{n}}$, $\widetilde J_n\subset \widetilde F_{l_{n}}$ such that
\roster
\item"(i)"
$F_{k_{n}} \widetilde J_n \subset \widetilde F_{l_{n}}$,
\item"(ii)" 
the mapping $F_{k_{n}}\times \widetilde J_n \ni (f, \widetilde f) \mapsto f\widetilde f \in \widetilde F_{l_{n}}$ is one-to-one, 
%(i.e. $(F_{k_{n+1}} - F_{k_{n+1}}) \cap (\widetilde J_n - \widetilde J_n) = \{0\}$),
\item"(iii)" 
%$\dfrac{\# \widetilde F_{l_{n}} - \# F_{k_{n}} \# \widetilde J_n}{\# \widetilde F_{l_{n}}} < 2\epsilon_n$,
%\item"(iv)" 
$\dfrac{\#\big((\widetilde J_n J_{n+1}) \Delta (C_{k_{n} + 1}\cdots C_{ k_{n+1}})\big)}{\# C_{k_{n} + 1}\cdots\#C_{ k_{n+1}}} < \dfrac1{2^n}$,
\item"(i)$'$" $\widetilde F_{l_{n}}  J_{n+1} \subset F_{k_{n+1}}$,
\item"(ii)$'$" 
the mapping $\widetilde F_{l_{n}}\times J_{n+1} \ni (\widetilde f, f) \mapsto \widetilde f f \in  F_{k_{n+1}}$ is one-to-one,
%(i.e. $(\widetilde F_{l_{n+1}} - \widetilde F_{l_{n+1}}) \cap (J_{n+1} - J_{n+1}) = \{0\}$),
\item"(iii)$'$" 
%$\dfrac{\# F_{k_{n+1}} - \# \widetilde F_{l_{n}} \# J_{n+1}}{\# F_{k_{n+1}}} < 2\epsilon_n$,
%\item"(iv)$'$" 
$\dfrac{\#\big((J_{n+1}  \widetilde J_{n+1}) \Delta (\widetilde C_{l_n + 1}\cdots \widetilde C_{ l_{n + 1}})\big)}{\# \widetilde C_{l_n + 1}\cdots\#\widetilde C_{ l_{n+1}}} < \dfrac1{2^n}$
\endroster
for each $n\ge 0$.
\endproclaim

Moreover, we show that each isomorphism intertwining $T$ with  $\widetilde T$ is a composition of seven ``elementary'' isomorphisms between  $(C,F)$-systems.
If $G=\Bbb Z$, $T$ and $\widetilde T$ are finite measure preserving, $F_n=\{0,1,\dots, h_n\}$, 
$\widetilde F_n=\{0,1,\dots, \widetilde h_n\}$ for some positive integers $h_n$ and $\widetilde h_n$
and every $n>0$ then Theorem~A provides a solution to the  Foreman--Rudolph--Weiss problem.

We also note that Theorem~A is a measure theoretical analogue of the classification of continuous $(C,F)$-actions on locally compact spaces obtained in  \cite{Da3, Theorem~A}.

By a {\it factor} of a probability preserving action we mean an invariant sub-$\sigma$-algebra of measurable subsets as well as the restriction of this action to this sub-$\sigma$-algebra.
It is well known that each factor of a rank-one transformation is of rank one. 
This  is no longer true for rank-one actions of general amenable groups (see \cite{DaVi} for counterexamples).
The following result is a description of all rank-one factors of a rank-one system with a finite invariant measure.
%However,  we can describe  all  rank-one actions that appear as factors of a given rank-one  action
% in the following way.

\proclaim{Theorem B}
Let $G$ be amenable.
Let $T$ and $\widetilde T$ be two finite measure preserving  $G$-actions associated with
 $(C_n,F_{n-1})_{n=1}^\infty$ and 
$(\widetilde C_n,\widetilde F_{n-1})_{n=1}^\infty$ respectively.
Then $\widetilde T$ is   isomorphic to a measure theoretical factor of $T$ if and only if there exist two increasing sequences
$
0=k_0<k_1<k_2<\cdots
$
and 
$
0=l_0<l_1<l_2<\cdots
$
 of non-negative integers and subsets $J_n\subset F_{k_{n}}$ such that
\roster
\item"(i)"
$\widetilde F_{l_n}  J_n \subset  F_{k_{n}}$,
\item"(ii)" 
the mapping $\widetilde F_{l_n}\times J_n \ni (\widetilde f,f) \mapsto \widetilde ff \in F_{k_{n}}$ is one-to-one, 
%(i.e. $(F_{k_{n+1}} - F_{k_{n+1}}) \cap (\widetilde J_n - \widetilde J_n) = \{0\}$),
\item"(iii)" 
$\dfrac{\# F_{k_n} - \# \widetilde F_{l_n} \#  J_n}{\#  F_{k_n}} < \frac 1{2^n}$,
\item"(iv)" 
$\dfrac{\#((J_{n-1}  C_{k_{n-1}+1}C_{k_{n-1}+2}\cdots C_{k_n}) \Delta \widetilde C_{l_{n-1}+1}\widetilde C_{k_{n-1}+2}\cdots \widetilde C_{l_n}J_{n})}
{\# \widetilde C_{l_{n-1}+1}\cdots \#\widetilde C_{l_n}\# J_{n}} < \frac 1{2^n}$
\endroster
for each $n\ge 1$.
\endproclaim

It is well known that the odometer $\Bbb Z$-actions are of rank one  \cite{dJ1}.
A description of the odometer factors for rank-one transformations was obtained recently in
\cite{Fo--We}  (see also  \cite{DaVi}).
We show how to deduce  this description  from Theorem~B.

A topological counterpart of Theorem~B is also proved: in Theorem~3.3, for an arbitrary discrete countable group $G$, we describe all proper continuous  factors of continuous  $(C,F)$-actions   of $G$ defined on locally compact Cantor spaces.

Denote by $\goth R_1^{\text{fin}}$ and $\goth R_1^{\infty}$ the spaces of parameters of the $(C,F)$-actions of $G$ with finite and infinite  invariant measure  respectively.
The two spaces have natural Polish topologies.
Define the  {\it isomorphism equivalence relation} on $\goth R_1^{\text{fin}}$ (and separately on $\goth R_1^{\infty}$) by saying that two $(C,F)$-sequences are isomorphic if the  $(C,F)$-actions associated with them  are measure theoretically isomorphic. 

\proclaim{Theorem C}
The isomorphism equivalence relation on $\goth R_1^{\text{fin}}$ is a 
$G_\delta$-subset of  $\goth R_1^{\text{fin}}\times \goth R_1^{\text{fin}}$.
In a similar way, the isomorphism equivalence relation on $\goth R_1^{\infty}$ is a 
$G_\delta$-subset of  $\goth R_1^{\infty}\times \goth R_1^{\infty}$.
\endproclaim

The first claim of Theorem~C extends and refines \cite{FoRuWe, Theorem~51}, where it was proved 
that the isomorphism equivalence relation is Borel
on the set of classical rank-one $\Bbb Z$-actions.
The proof of \cite{FoRuWe, Theorem~51} is based heavily on  King's weak closure theorem
\cite{Ki}, which does not hold for actions of arbitrary amenable groups (see \cite{DoKw} for counterexamples in the case where $G=\Bbb Z^2$).
Our proof of Theorem~C, given for actions of either arbitrary amenable groups (the first claim)  or arbitrary countable   groups (the second claim), is based solely on Theorem~A.

The outline of our paper is the following.
In Section~1 we remind definitions related to rank-one actions and the $(C,F)$-construction.
The ``elementary'' isomorphisms between $(C,F)$-actions: calibration, telescoping, reduction and chain equivalence are explained in detail.
Theorem~A is proved in Section~2.
Section~3 describes  rank-one factors of rank-one probability preserving systems.
We first prove a topological counterpart of Theorem~B, then Theorem~B itself and finally show that Theorem~B implies the description of odometer factors from \cite{Fo--We} and
\cite{DaVi}.
Section~4 is devoted to the proof of Theorem~C.
\vskip 2 pt

{\sl Acknowledgements.}  The authors thank B. Weiss for his useful remarks.

\head 1. Rank-one actions and  $(C,F)$-construction
\endhead

\subhead 1.1. Group actions of rank one
\endsubhead
Let $G$ be a discrete infinite countable  group. 
Let $T=(T_g)_{g\in G}$ be a free measure preserving action of $G$ on a standard  non-atomic 
$\sigma$-finite measure space $(X,\goth B,\mu)$.
By a {\it  Rokhlin tower for $T$} we mean a pair $(B,F)$, where $B\in\goth B$ with $\mu(B)>0$ and $F$ is a finite subset of $G$ with $1_G\in F$ such that
the subsets $T_fB$, $f\in F$, are mutually disjoint.
We let $X_{B,F}:=\bigsqcup_{f\in F}T_fB\in\goth B$.
By $\xi_{B,F}$ we mean the finite partition of  $X_{B,F}$ into the subsets $T_fB$, $f\in F$.
If $x\in T_fB$ then we set $O_{B,F}(x):=\{T_gx\mid g\in Ff^{-1}\}$.

\definition{Definition 1.1 \cite{DaVi}} Let $\{1_G\}=F_0\subset F_1\subset F_2\subset\cdots$ be an increasing sequence
of finite subsets in $G$.
We say that $T$ is of {\it rank-one along $(F_n)_{n=0}^\infty$} if there is a decreasing sequence 
$B_0\supset B_1\supset\cdots$ of subsets of positive measure in $X$ such that
$(B_n,F_n)$ is a
 Rokhlin tower for $T$ for each $n\in\Bbb N$ and
\roster
\item"(i)" $\xi_{B_n,F_n}\prec\xi_{B_{n+1},F_{n+1}}$ for each $n\ge 0$ and $\bigvee_{n=0}^\infty \xi_{B_n,F_n}$ is the partition of $X$ into singletons (mod 0),
\item"(ii)" $\{T_gx\mid g\in G\}=\bigcup_{n=1}^\infty O_{B_n,F_n}(x)$ for a.e. $x\in X$.
\endroster
\enddefinition

It follows from (i) that $X_{B_0,F_0}\subset X_{B_1,F_1}\subset X_{B_2,F_2}\subset\cdots$ and $\bigcup_{n=0}^\infty X_{B_n,F_n}=X$.

In the case where $G=\Bbb Z$, a classical rank-one transformation  corresponds to the rank-one along
a sequence $(\{0,1,\dots, h_n-1\})_{n=1}^\infty$ for some increasing sequence of integers $h_1<h_2<\cdots$
according to Definition~1.1.
The rank-one $\Bbb Z$-actions along  arbitrary sequences $(F_n)_{n=0}^\infty$ correspond to the class of funny rank-one transformations (see \cite{Fe} for the finite measure preserving case).
We also note that in the classical case of finite measure preserving rank-one transformations, ~(ii)
follows from (i).
However, in the general case, it cannot be omitted: see a counterexample  \cite{DaVi, Example~4.4} for the free group with 2 generators.

\subhead 1.2. $(C,F)$-spaces, tail equivalence relations and  return time cocycles 
\endsubhead
We remind the $(C,F)$-construction as it appeared in \cite{Da2}.
Let $G$ be  a discrete countable group.
Let $\Cal T=(C_n,F_{n-1})_{n=1}^\infty$ be  a sequence of (pairs of) finite subsets  of $G$ such that
$\#F_{0} = 1$ and for 
 each $n>0$,
$$
\aligned
	  &  %1_G\in F_n\cap C_n, \  
	  \#C_{n} > 1, \\ 
      &F_{n} C_{n+1}\subset F_{n+1},\\ 
      &F_{n} c\cap F_{n} c' = \emptyset\text{ if $c, c'\in C_{n+1}$ and $c \neq c'$. }
\endaligned
\tag1-1
$$
We let $X_n := F_{n} \times C_{n+1} \times C_{n+2} \times\ldots$ and endow this set with the infinite product topology. 
Then $X_n$ is a compact Cantor space. The mapping 
$$
	X_n \ni (f_n,c_{n+1},c_{n+2}\ldots) \mapsto (f_n c_{n+1},c_{n+2},\ldots) \in X_{n+1}\tag1-2
$$
is a continuous embedding of $X_n$ into $X_{n+1}$. 
Therefore the topological  inductive limit $X$ of the sequence $(X_n)_{n\geq 0}$  is  well defined.
Moreover, $X$ is a locally compact Cantor space. 
We call $X$ {\it the $(C,F)$-space associated with $\Cal T$}.
It is convenient to consider $X$ as the union $\bigcup_{n=1}^\infty X_n$ of the increasing
sequence $X_0\subset X_1\subset\cdots$ of compact open subsets, where the corresponding embeddings are given by \thetag{1-2}.
For a subset $A\subset F_n$, we let
$$
	[A]_n := \{x=(f_n,c_{n+1},\ldots)\in X_n, f_n\in A\}
$$
and call this set an $n$-{\it cylinder} in $X$.
It is open and compact in $X$. 
Every open subset of $X$ is a union of cylinders.
For brevity, we will write $[f]_n$ for $[\{f\}]_n$ for an  element $f\in F_n$.

Two points $x=(f_n,c_{n+1},\ldots)$ and $x'=(f_n',c_{n+1}',\ldots)$ of $X_n$ are called {\it tail equivalent}
if there is $N>n$ such that $c_l=c_{l}'$ for  each $l>N$.
We thus obtain the tail equivalence relation on $X_n$.
The {\it  tail equivalence relation}  $\Cal R$ on $X$ is defined as follows: for each $n\ge 0$, the restriction of $\Cal R$ to $X_n$ is 
the tail equivalence relation on $X_n$. 
We note that $\Cal R$ is {\it Radon uniquely ergodic}, i.e. there is a unique  $\Cal R$-invariant Radon measure $\mu$ on $X$ such that $\mu(X_0)=1$.
We call it  the {\it  Haar measure  for $\Cal R$}.
It is $\sigma$-finite.
Let $\kappa_n$ be the equidistribution on $C_n$.
We define a measure $\nu_n$ on $F_n$ by setting $\nu_0(F_0)=1$ and
 $\nu_n(\{f\})=1/\prod_{k=1}^n \#C_k$ for each  $f\in F_n$ and $n> 0$.
 Then
 $$
 \mu\restriction X_n=\nu_n\otimes\bigotimes_{k>n}\kappa_k\quad\text{for each $n>0$}.
 $$
The Haar measure for $\Cal R$ is finite if and only if 
$$
\prod_{n=1}^\infty\frac{\# F_{n+1}}{\# F_n\# C_{n+1}}<\infty.\tag1-3
$$
Of course, $\Cal R$ is {\it minimal} on $X$, i.e. each $\Cal R$-class is dense in $X$.

Define a cocycle $\alpha:\Cal R\to G$ by setting
$$
\alpha(x,\widetilde x):=\lim_{m\to\infty}f_nc_{n+1}\cdots c_m\widetilde c_m^{-1}\cdots\widetilde c_{n+1}^{-1}\widetilde f_n^{-1}
$$
if $x=(f_n,c_{n+1},c_{n+2},\dots)\in X_n$ and 
$\widetilde  x=(\widetilde f_n,\widetilde c_{n+1},\widetilde c_{n+2},\dots)\in X_n$
for some $n\ge 0$.
It is straightforward to verify that  $\alpha$ is well defined and it satisfies the cocycle identity
$$
\alpha(x,\widetilde x)\alpha(\widetilde x,\widehat x)=\alpha(x,\widehat x)
$$ 
for all $x,\widetilde x,\widehat x\in X$
such that $(x,\widetilde x)\in\Cal R$ and $(\widetilde x,\widehat x)\in\Cal R$.
We call $\alpha$ {\it the return time cocycle} of $\Cal R$.

\subhead 1.3. $(C,F)$-actions: topological and measure theoretical
\endsubhead
Given $g\in G$,  let 
$$
X_n^g:=\{(f_n, c_{n+1}, c_{n+2},\dots)\in X_n\mid gf_n\in F_n\}.
$$
Then $X_n^g$ is a compact open subset of $X_n$ and $X_n^g\subset X^g_{n+1}$.
  Hence the union
  $X^g:=\bigcup_{n\ge 0}X_n^g$ is an open subset of $X$.
Let $X^G:=\bigcap_{g\in G}X^g$.
 Then $X^G$ is a $G_\delta$-subset of $X$.
  Hence $X^G$ is Polish and totally disconnected  in the induced topology. 
  Given $g\in G$ and  $x\in X_G$, there is $n >0$ such that
  $x= (f_n, c_{n+1},\dots )\in X_n$ and $gf_n\in F_n$. 
  We now let 
  $$
  T_gx:= (gf_n, c_{n+1}, \dots)\in X_n\subset X.
  $$
  It is straightforward  to verify that
  \roster
  \item"---" $T_gx\in X^G$,
  \item"---" the mapping $T_g:X^G\ni x\mapsto T_gx\in X^G$ is a homeomorphism of $X^G$ and
  \item"---" $T_gT_{g'}=T_{gg'}$ for all $g, g'\in G$,
  \item"---"  $\alpha(T_gx,x)=g$ for all $g\in G$ and $x\in X^G$,
  where $\alpha$ is the return time cocycle of $\Cal R$.
\endroster
Hence, $T:=(T_g)_{g\in G}$
is a continuous, well defined $G$-action on $X^G$.

\definition{Definition 1.2 \cite{Da2}} The action $T$ is called {\it the topological $(C,F)$-action of $G$ associated with $\Cal T$}.
\enddefinition

The topological $(C,F)$-action is free.
The subset $X^G$ is $\Cal R$-invariant.
The $T$-orbit equivalence relation 
coincides with the restriction of $\Cal R$ to $X^G$.
It was shown in \cite{Da2} that
$X^G=X$ if and only if for each $g\in G$ and $n >0$, there is
$m > n$ such that
$$
gF_nC_{n+1}C_{n+2}\cdots C_m\subset F_m.
\tag1-4
$$
Thus, if \thetag{1-4} holds then $T$ is a minimal continuous $G$-action on a locally compact Cantor space $X$.
Moreover, $T$ is {\it Radon uniquely ergodic}, i.e. there exists a unique  $T$-invariant Radon measure $\xi$ on $X$ such that $\xi(X_0)=1$.
Of, course $\xi$ is the Haar measure for  $\Cal R$.
If  $T$ is Radon uniquely ergodic and \thetag{1-3} holds then $T$ is uniquely ergodic in the classical sense.

From now on, $T$ is a topological $(C,F)$-action of $G$ on $X^G$ and
 $\mu$ is the Haar measure for $\Cal R$. 
 Since $X^G$ is $\Cal R$-invariant, we obtain that either $\mu(X^G)=0$ or $\mu(X\setminus X^G)=0$.
In the latter case $T$ is  conservative and ergodic.
The following two results were obtained  in \cite{Da2}:

\proclaim{Fact A}
$\mu(X\setminus X^G)=0$ if and only if for each $g\in G$ and every $n \ge0$,
$$
\lim_{m\to\infty}\nu_m\big((gF_nC_{n+1}C_{n+2}\cdots C_m)\cap F_m\big)=\nu_n(F_n).\tag1-5
$$
\endproclaim

\proclaim{Fact B}
If  $\mu(X\setminus X^G)=0$ and $\mu(X)<\infty$  then
$G$ is amenable and 
$(F_n)_{n=1}^\infty$ is a left F{\o}lner sequence in $G$.
\endproclaim

We also note that if $G$ admits a finite measure preserving $(C,F)$-action then
the F{\o}lner sequence $(F_n)_{n=1}^\infty$ possesses the following ``near tiling'' property:  for each pair  of integers $m>n>0$, there is a finite subset $D_{n,m}$ such that $F_nD_{n,m}\subset F_m$, $F_nd\cap F_nd'=\emptyset$ for all $d\ne d'\in D_{n,m}$ and for each $\epsilon>0$, if $n$ is large enough then
$\#(F_nD_{n,m})/\#F_m>1-\epsilon$ for every $m>n$.
We do not know whether each countable amenable group has a F{\o}lner sequence with the near tiling property.

\definition{Definition 1.3} If $\mu(X\setminus X^G)=0$ then the dynamical system $(X,\mu, T)$ (or simply $T$) is called
{\it the  measure preserving $(C,F)$-action associated with  $\Cal T$.}
\enddefinition

\comment

If $T$ is a measure preserving $(C,F)$-action then each  $T$-invariant probability measure $\xi$ on $X^G$
is also a probability measure on $X$ which is invariant under the tail equivalence relation $\Cal R$ on $X$.
Since $\Cal R$ is Radon uniquely ergodic %and \thetag{1-3} holds,
 we deduce that
$\xi$ is a scaling of the   Haar measure for $\Cal R$.
Thus, $T$ is uniquely ergodic.

\endcomment

\proclaim{Fact C (\cite{Da2}, \cite{DaVi})} Each measure preserving $(C,F)$-action  is of rank one.
Conversely, each rank-one  $G$-action is  isomorphic (via a   measure preserving isomorphism) to a $(C,F)$-action.
\endproclaim

We illustrate this isomorphism with the example of classical $\Bbb Z$-actions of rank one.

\example{Example 1.4} Let $(r_n)_{n=1}^\infty$ be a sequence of natural numbers, $r_n>1$ for each $n\in\Bbb N$.
Let $s_n:\{0,1,\dots,r_n-1\}\to\Bbb Z_+$ be a sequence of mappings.
Then there is a geometric cutting-and-stacking inductive construction with a single tower on each step to craft a measure preserving rank-one transformation $S$ of an interval $[0,\alpha)$ furnished with the Lebesgue measure.
On the $n$-th step of the construction, we have an $n$-th tower consisting of $h_n$ levels numbered by $0,1,\dots,h_n-1$ from bottom to  top.
Every level is a semi-interval $[a,b)\subset [0,+\infty)$ of length $1/(r_1\cdots r_n)$.
The rank-one transformation $S$ is defined partially on the $n$-th tower: it moves each level of the tower, except the highest one,  up
one level.
 We cut the $n$-th tower  into $r_n$ subtowers numbered with $0,1,\dots, r_n-1$ from the left to the right.
 Thus, each subtower consists of semi-intervals of length 
 $1/(r_1\cdots r_{n+1})$.
Then $s_n(i)$ ``spacers''  are put on the top of the $i$-th subtower, $i=0,\dots, r_n-1$.
Each spacer is also a semiinterval
of length 
 $1/(r_1\cdots r_{n+1})$.
We note that  $(r_k)_{k=1}^\infty$ is called {\it the sequence of cuts} and $(s_k)_{k=1}^\infty$ is called {\it the sequence of spacer maps}.
Then we stack the subtowers (extended with the spacers) into a single $(n+1)$-th tower by putting the $i$-th
one on the top of the $(i-1)$-th one, $i=1,\dots, r_n-1$.
If the $n$-th tower is of height $h_n$ then the height  of the $(n+1)$-th tower
equals
$$
h_{n+1}=r_nh_n+\sum_{i=0}^{r_n-1}s_n(i).
$$
We number the levels in the $(n+1)$-th tower with $0,1,\dots,h_{n+1}-1$ from bottom to top and define $S$
partially on the $(n+1)$-th tower according to this numbering.
And so on.
``At the end'' of the inductive construction, $S$ is determined almost everywhere on $[0,\alpha)$.
We also note that $\alpha=\lim_{n\to\infty}\frac {h_n}{r_1\cdots r_n}$.
For details of this cutting-and-stacking construction we refer to  \cite{Si}.
We now define two sequences of finite subsets in $\Bbb Z$ by setting $F_0:=\{0\}$ and for each $n>0$,
$$
\align
F_n&:=\{0,1,\dots,h_n-1\} \quad\text{and}\\
C_{n+1}&:=\bigg\{0, h_n+s_n(0), 2h_n+s_n(0)+s_n(1), \dots,(r_n-1)h_n+\sum_{i=0}^{r_n-2}s_n(i)\bigg\}.
\endalign
$$
Then the sequence $\Cal T:=(F_n,C_{n+1})_{n=0}^\infty$ satisfies \thetag{1-1}
and \thetag{1-5}.
Moreover, $\alpha<\infty$ if and only if \thetag{1-3} is satisfied.
Hence,  the $(C,F)$-action $T=(T_m)_{m\in\Bbb Z}$ of $\Bbb Z$ is well defined.
It is straightforward to verify  $S$ is isomorphic to $T_1$.
% according to~Fact~C.
\endexample

\subhead 1.4. Calibrations\endsubhead
If $\Cal T$ satisfies \thetag{1-1} %and  \thetag{1-3} 
and $\boldsymbol z:=(z_n)_{n\ge 1}$ is an arbitrary sequence of elements of $G$, we let
$$
C_n':=z_n^{-1}C_nz_{n+1}\quad\text{and}\quad F_{n-1}':=F_{n-1}z_{n}\quad\text{ for each $n\ge 1$}
$$
 and $\Cal T':=(C_n',F_{n-1}')_{n=1}^\infty$.
Then $\Cal T'$ also satisfies \thetag{1-1}. %and \thetag{1-3}.
We call $\Cal T'$  the  {\it $\boldsymbol z$-calibration} of $\Cal T$.
Let $X$ and  $X'$ be the $(C,F)$-spaces associated with $\Cal T$ and $\Cal T'$ respectively.
Denote by $\Cal R$ and ${ \Cal R}'$ the tail equivalence relations on $X$ and $ X'$
respectively.
We define  a mapping $\phi_{\boldsymbol z}:X\to X'$ by setting
$$
\phi_{\boldsymbol z}(f_n,c_{n+1},c_{n+2},\dots)=(f_nz_{n+1}, z_{n+1}^{-1}c_{n+1}z_{n+2}, z_{n+2}^{-1}c_{n+2}z_{n+3},\dots )\in X_n'\subset X'
$$
whenever $(f_n,c_{n+1},c_{n+2},\dots)\in X_n\subset X $, for each $n\ge 0$.
It is straightforward to verify that  $\phi_{\boldsymbol z}$ is a  homeomorphism.
Moreover, $\phi_{\boldsymbol z}$ maps bijectively each $\Cal R$-class in $X$ onto
an $\Cal R'$-class in $ X'$ and
$$
\alpha'(\phi_{\boldsymbol z}x, \phi_{\boldsymbol z}x')=\alpha(x,x')\quad\text{for each $(x,y)\in\Cal R$},\tag1-6
$$
where $\alpha$ and $\alpha'$ denote the return time cocycles of $\Cal R$ 
and $\Cal R'$ respectively.
We  call $\phi_{\boldsymbol z}$  the  {\it $\boldsymbol z$-calibration mapping}.
Of course, $\phi_{\boldsymbol z}$ maps the Haar measure on $X$ onto the Haar measure on $X'$.

Choosing $\boldsymbol z$ in an appropriate way we may  assume without loss of 
generality\footnote{This means that we can modify (calibrate) $\Cal T$ so that  the $(C, F) $-action associated 
with the modified sequence is isomorphic to the original $(C,F)$-action.}
 that the  condition
$$
1_G\in\bigcap_{n\ge 0}F_n\cap\bigcap_{n\ge 1}C_n
\tag 1-7
$$
is satisfied for $\Cal T$ in addition to \thetag{1-1}. 

If $\Cal T$ satisfies  \thetag{1-4} or \thetag{1-5} then 
$\Cal T'$ also satisfies~\thetag{1-4} or \thetag{1-5} respectively.
Hence,  the $(C,F)$-actions $T$ and $T'$ associated with $\Cal T$ and $\Cal T'$ respectively
are well defined simultaneously.
It follows from  \thetag{1-6} that $T$ and $T'$ are conjugate via $\phi_{\boldsymbol z}$, i.e.
$\phi_{\boldsymbol z}T_g\phi_{\boldsymbol z}^{-1}=T_g'$ for each $g\in G$.

\subhead  1.5. Telescopings 
\endsubhead
Let a sequence $\Cal T=(C_n,F_{n-1})_{n=1}^\infty$ satisfy   \thetag{1-1}. % and \thetag{1-3}.
Given a strictly increasing  infinite sequence of integers $\boldsymbol l=(l_n)_{n=0}^\infty$ such that
$l_0=0$, we let 
$$
\widetilde F_n:=F_{l_n},\quad \widetilde C_{n+1}:=C_{l_n+1}\cdots C_{l_{n+1}}
$$
for each $n\ge 0$.
The sequence $\widetilde{\Cal T}:=(\widetilde C_n,\widetilde F_{n-1})_{n=1}^\infty$ is called the {\it $\boldsymbol l$-tele\-scoping} of $\Cal T$ \cite{Da2}.
It is easy to check that $\widetilde{\Cal T}$   satisfies  \thetag{1-1}.  %and  \thetag{1-3}.
Let $X$ and $\widetilde X$ denote the $(C,F)$-spaces associated with $\Cal T$
and $\widetilde{\Cal T}$ respectively.
Denote by $\Cal R$ and $\widetilde{\Cal R}$ the tail equivalence relations on $X$ and $\widetilde X$
respectively.
There is a canonical mapping $\iota_{\boldsymbol l}$ of  $X$ onto $\widetilde X$ associated with $\boldsymbol l$.
If $x\in X$ then we select the smallest  $n\ge 0$ such that $x=(f_{l_n},c_{l_n+1},c_{l_n+1},\dots)\in X_{l_n}$ and put
$$
\iota_{\boldsymbol l}(x):=(f_{l_n},c_{l_n+1}\cdots c_{l_{n+1}}, c_{l_{n+1}+1}\cdots c_{l_{n+2}},\dots)\in \widetilde X_n\subset\widetilde X,
$$
where $\widetilde X_n=\widetilde F_n\times\widetilde C_{n+1}\times \widetilde C_{n+2}\times\cdots$.
It is routine to verify that 
\roster
\item"---" $\iota_{\boldsymbol l}$ is a homeomorphism of $X$
onto $\widetilde X$,
\item"---" $\iota_{\boldsymbol l}$ maps bijectively each $\Cal R$-class in $X$ onto
an $\widetilde{ \Cal R}$-class in $X$,
\item"---" 
$\iota_{\boldsymbol l}$ transfers the Haar measure for $\Cal R$ to the Haar measure
for $\widetilde {\Cal R}$.
\endroster
Moreover,
$$
\widetilde\alpha(\iota_{\boldsymbol l}x,\iota_{\boldsymbol l}x')=\alpha(x,x')\quad\text{for each $(x,x')\in\Cal R$},\tag{1-8}
$$
where $\alpha$ and $\widetilde\alpha$ are the return time cocycles of $\Cal R$ and $\widetilde{\Cal R}$ respectively.
We call $\iota_{\boldsymbol l}$ {\it the $\boldsymbol l$-telescoping mapping.}

If $\Cal T$ satisfies  \thetag{1-4} or \thetag{1-5} then 
$\widetilde{\Cal T}$ also satisfies~\thetag{1-4} or \thetag{1-5} respectively.
Hence,  the $(C,F)$-actions $T$ and $\widetilde T$ associated with $\Cal T$ and $\widetilde{\Cal T}$ respectively
are well defined.
It follows from  \thetag{1-8} that $T$ and $\widetilde T$ are conjugate via $\iota_{\boldsymbol l}$, i.e.
$\iota_{\boldsymbol l}T_g\iota_{\boldsymbol l}^{-1}=\widetilde T_g$ for each $g\in G$.

\subhead 1.6. Reductions
\endsubhead
Let a sequence $\Cal T=(C_n,F_{n-1})_{n=1}^\infty$ satisfy   \thetag{1-1}.
% and \thetag{1-3}.
Let $\boldsymbol A:=(A_n)_{n=1}^\infty$ be  a sequence of nonempty subsets $A_n\subset C_{n}$ such that
%$1_G\in A_n$ for each $n\in\Bbb N$ and 
$$
\sum_{n=1}^\infty(1- \kappa_n(A_n))<\infty.
$$ 
We will assume that $A_n$ is a proper subset of $C_n$ for infinitely many $n$.
Denote by $ \kappa_n^*$
the equidistribution 
on $A_n$ for each $n\in\Bbb N$.
Let $\Cal T^*:=(A_n,F_{n-1})_{n=1}^\infty$.
The sequence $\Cal T^*$  is called   {\it the $\boldsymbol A$-reduction} of $\Cal T$
\cite{Da2}.
It is easy to check that $\Cal T^*$   satisfies  \thetag{1-1}.
% and~\thetag{1-3}.
%It was shown in \cite{Da2} that   $\Cal T'$   satisfies \thetag{1-5}.
Let $X$ and $X^*$ be the $(C,F)$-spaces associated with $\Cal T$ and $\Cal T^*$.
Denote by 
$\Cal R$ and $\Cal R^*$ the tail equivalence relations on $X$ and $ X^*$
respectively.
Let $\mu$ and $\mu^*$ denote the Haar measures on $X$ and $X^*$
respectively.
We note that for each $n\ge0$, the identity mapping embeds the set 
$$
X_n^*:=F_n\times A_{n+1}\times A_{n+2}\times\cdots
$$
into  $X_n:=F_n\times C_{n+1}\times C_{n+2}\times\cdots$.
Hence, we can consider $X_n^*$ 
as a nowhere dense closed subset of $X_n$.
It follows that $X^*=\bigcup_{n\ge 0}X_n^*$ embeds naturally  into $X$ as an $F_\sigma$-subset of the first Baire category.
Of course, $X^*$ is $\Cal R$-invariant and, hence, dense in $X$.
The restriction of $\Cal R$ to $X^*$ is $\Cal R^*$.
We note that
$$
\mu(X^*)\ge \mu(X_n^*)>\prod_{j>n}{\kappa_j(A_j)}>0.
$$
Since $X^*$ is $\Cal R$-invariant and $\mu$ is $\Cal R$-ergodic, it follows that $\mu(X\setminus X^*)=0$.
Thus, $X^*$ is of full measure in $X$.
There is a  canonical measure scaling Borel isomorphism $\rho_{\boldsymbol A}$ of  $(X,\mu)$ onto $(X^*,\mu^*)$:
$$
\rho_{\boldsymbol A}x:=x\quad\text{if $x\in X^*$}.
$$
The reader should not confuse $x$ from the lefthand side ($x$ is a point of the $(C,F)$-space $X$) with $x$ from the righthand side ($x$ is a point of the $(C,F)$-space $X^*$).
Thus, $\rho_{\boldsymbol A}$ is defined only on $X^*$, which is a
$\mu$-conull  $F_\sigma$- subset of  $X$.
It is straightforward to verify that
\roster
\item"---" the inverse mapping $\rho_{\boldsymbol A}^{-1}:X^*\to X$ 
is continuous, 
\item"---"
$\rho_{\boldsymbol A}$ maps bijectively each $\Cal R$-class in $X^*\subset X$ onto
an $\Cal R^*$-class in $X^*$,
\item"---"
$
\frac{d(\mu\circ \rho_{\boldsymbol A}^{-1})}{d\mu^*}=\prod_{m>0}{\kappa_m(A_m)}
$
almost everywhere
and
\item"---"  
$\alpha^*(\rho_{\boldsymbol A}x,\rho_{\boldsymbol A}x')=\alpha(x,x')$ for each $(x,x')\in\Cal R\cap(X^*\times X^*)$,
\endroster
where $\alpha$ and $\alpha^*$ are the return time cocycles of $\Cal R$ and $\Cal R^*$ respectively.
We call $\rho_{\boldsymbol A}$ {\it the $\boldsymbol A$-reduction mapping.}

If $\Cal T$ satisfies  \thetag{1-4} or \thetag{1-5} then 
$\Cal T^*$ also satisfies~\thetag{1-4} or \thetag{1-5} respectively.
Hence,  the $(C,F)$-actions $T$ and $ T^*$ associated with $\Cal T$ and $\Cal T^*$ respectively
are well defined.
Moreover,  $T$ and $T^*$ are conjugate via $\rho_{\boldsymbol A}$, i.e.
$\rho_{\boldsymbol A}T_g\rho_{\boldsymbol A}^{-1}= T_g^*$ a.e. for each $g\in G$.

\comment

Indeed,
fix $n\in\Bbb N$.
Since $\mu\restriction X_n=\nu_n\otimes\kappa_{n+1}\otimes\kappa_{n+2}\otimes\cdots$, it follows from the Borel-Cantelli lemma that for a.e. $x=(f_n,c_{n+1},c_{n+2},\dots)\in X_n$, there is
$N=N_x>n$ such that $c_m\in A_m$ for each $m>N$.
We then let
$$
\iota_{\boldsymbol A,n}(x):=(f_nc_{n+1}\cdots c_N, c_{N+1},c_{N+2},\dots)\in X_N'\subset X'.
$$
It is routine to verify that $\iota_{\boldsymbol A,n}:X_n\ni x\mapsto \iota_{\boldsymbol A,n}(x)\in X'$ is a well defined 
 mapping.
 Moreover,
$\iota_{\boldsymbol A,n+1}\restriction X_n=\iota_{\boldsymbol A,n}$ for each $n\in\Bbb N$.
Hence a measurable mapping $\iota_{\boldsymbol A}:X\to X'$ is a well defined
by the restrictions $\iota_{\boldsymbol A}\restriction X_n=\iota_{\boldsymbol A,n}$ for all $n\in\Bbb N$.
It is straightforward to verify that $\iota_{\boldsymbol A}$ is a measure scaling  isomorphism of $(X,\mu)$ onto $(X',\mu')$ with  
$\mu\circ \iota_{\boldsymbol A}^{-1}=\big(\prod_{m>0}{\kappa_m(A_m)}\big)\cdot\mu'$ and
 $\iota_{\boldsymbol A}T_g=T_g'\iota_{\boldsymbol A}$ for each $g\in G$.
\endcomment

\proclaim{Fact D \cite{Da2}} Let $\Cal T$ be a $(C,F)$-sequence satisfying \thetag{1-1}, \thetag{1-3} and \thetag{1-5}.
Then there is a telescoping $\widetilde {\Cal T}$ of $\Cal T$ 
and a reduction  $\widetilde {\Cal T}^*$ of  \, $\widetilde {\Cal T}$ such that 
$\widetilde {\Cal T}^*$ satisfies~\thetag{1-1}, \thetag{1-3} and \thetag{1-4}.
\endproclaim

It follows from Facts~C and D that  each finite measure preserving rank-one action of $G$
 is measure theoretically isomorphic 
to a minimal uniquely ergodic continuous $(C,F)$-action on a locally compact Cantor space.

\subhead 1.7. Chain equivalence 
\endsubhead
%Let $T$ and $T'$ be two $(C,F)$-actions associated with sequences 
The chain equivalence for $(C,F)$-systems was introduced implicitly (without any name) in the proof of \cite{Da3, Theorem~A}.
We consider here a slightly more general version of that concept.
Let $\Cal T=(C_n,F_{n-1})_{n=1}^\infty$ and $\Cal T'=(C_n',F_{n-1}')_{n=1}^\infty$ be two $(C,F)$-sequences satisfying
\thetag{1-1}.
% and \thetag{1-3}.
Denote by $X$ and $X'$ the corresponding $(C,F)$-spaces.

\definition{Definition 1.5} 
We say that $\Cal T$  is {\it chain equivalent} to  $\Cal T'$ if there exist sequences
$\boldsymbol A:=(A_n)_{n=0}^\infty$ and $\boldsymbol B:=(B_n)_{n=1}^\infty$
of finite subsets in $G$
such that for each $n\ge 1$,
$$
\gather
A_{n-1}B_n=C_n,\quad B_nA_{n}=C'_{n},\tag1-9\\
F'_{n-1}B_n\subset F_n,\quad  F_{n-1}A_{n-1}\subset F'_{n-1}\quad\text{and}\tag1-10
\\
A_{n-1}^{-1}A_{n-1}\cap B_nB_n^{-1}=B_n^{-1}B_n\cap
A_{n}A_{n}^{-1}=\{1_G\}.
\tag1-11
\endgather
$$
\enddefinition

\comment %%%%%%%%%%%%%%%%%%%%%%%%%%%%%%%%

\proclaim{Lemma 1.6} 
Let $\Cal T$ and $\Cal T'$ be congruent. 
Then there exist  calibrations of $\Cal T$ and $\Cal T'$ which are congruent and for which
the corresponding sequences  $\widetilde{\boldsymbol A}=(\widetilde A_n)_{n=1}^\infty$ and 
$\widetilde{\boldsymbol B}=(\widetilde B_n)_{n=1}^\infty$ (for which~\thetag{1-9}--\thetag{1-11} hold)
satisfy an additional property
$$
1\in \widetilde A_n\cap \widetilde B_n\quad\text{for each $n\in\Bbb N$.}\tag1-12
$$
Hence, \thetag{1-7} holds for these calibrations of $\Cal T$ and $\Cal T'$.
\endproclaim

\demo{Proof}
We will choose successively $z_1,v_1,z_2,v_2,\dots$ such that
$$
F_0z_1=\{1_G\}, \ z_1^{-1}A_1v_1\ni 1_G, \   v_1^{-1}B_1z_2\ni 1_G, \  z_2^{-1}A_2v_2\ni 1_G,\ 
 v_2^{-1}B_2z_3\ni 1_G,\dots.
$$
Let $\boldsymbol z:=(z_n)_{n=1}^\infty$
and $\boldsymbol v:=(v_n)_{n=1}^\infty$.
Denote by $\widetilde{\Cal T}$ and $\widetilde{\Cal T'}$ the $\boldsymbol z$-calibration of $\Cal T$  and the
$\boldsymbol v$-calibration of $\Cal T'$ respectively.
Then  $\widetilde{\Cal T}=(z_n^{-1}C_nz_{n+1},F_{n-1}z_n)_{n=1}^\infty$ and 
$\widetilde{\Cal T'}=(v_n^{-1}C_n'v_{n+1},F_{n-1}'v_n)_{n=1}^\infty$.
We let for each $n\in\Bbb N$,
$$
\widetilde A_n:=z_n^{-1}A_nv_n\quad\text{and}\quad \widetilde B_n:=v_n^{-1}B_n
z_{n+1}.
$$
It is a routine to verify that \thetag{1-9}--\thetag{1-11} hold for $\widetilde A_n$, 
$\widetilde B_n$, $z_n^{-1}C_nz_{n+1}$, $F_{n-1}z_n$, $v_n^{-1}C_n'v_{n+1}$, $F_{n-1}'v_n$  in place of $A_n$, $B_n$, $C_n$, $F_{n-1}$, $C_n'$, $F_{n-1}'$ respectively.
This means that $\widetilde{\Cal T}$ and $\widetilde{\Cal T'}$ are congruent.
Of course, $1\in\widetilde A_n\cap\widetilde B_n$ for each $n\in\Bbb N$.
\qed
\enddemo

We will call the quadruple $(\boldsymbol  z, \boldsymbol v, \widetilde{\boldsymbol A}, \widetilde{\boldsymbol B})$ {\it a normalization for} $(\boldsymbol A,\boldsymbol B)$.

Suppose that $\Cal T$ and $\Cal T'$ are congruent and
\thetag{1-9}--\thetag{1-12} hold.
Let $\Cal R$ and $\Cal R'$ denote the tail equivalence relations on $X$ and $X'$
respectively.
Denote by $\alpha$ and $\alpha'$ the return time cocycles of  $\Cal R$ and $\Cal R'$ respectively.

We now define a mapping $\phi$ from  $[1_G]_0$ in $X$ to $[A_1]_1$ in $X'$.
Take $x\in[1_G]_0\subset X'$.
In view of the equality in the  lefthand side in \thetag{1-9} and \thetag{1-11}, there is a unique expansion
$x=(1_G,a_1b_1,a_2b_2,\dots )$ of $x$ such that $a_n\in A_n$ and $b_n\in B_n$ for each $n\in\Bbb N$.
Then we let
$$
\phi(x):=(a_1, b_1a_2,b_2a_3,\dots).
$$
%Since $A_1\subset F_0'=\{1_G\}$ and the equality in the righthand side in \thetag{1-7} holds,
 It follows  that $\phi(x)\in[A_1]_1\subset  X'$.
 We note that $\phi$ is one-to-one thanks to~\thetag{1-11}.
Of course, $\phi$ is continuous.
Hence $\phi$ is a homeomorphism of
$[1_G]_0$ in $X$ onto $[A_1]_1$ in $X'$.
It is straightforward to deduce from \thetag{1-9} and \thetag{1-11} that two points $x,\widetilde x\in [1_G]_0\subset X$ are tail equivalent in $X$ if and only if $\phi(x)$ and $\phi(\widetilde x)$ are tail equivalent in $X'$.
Moreover,
$$
\alpha(x,\widetilde x )=\alpha'(\phi(x),\phi(\widetilde x))
$$
for all $(x,\widetilde x)\in\Cal R$.

We now extend $\phi$ to a homeomorphism  $\psi_{\boldsymbol A,\boldsymbol B}$ from $X$ to $X'$.
Take $x\in X$.
Then there is $n\ge 0$ such that 
$x=(f_n,c_{n+1},c_{n+2}, \dots)\in X_n$.
We now let  
$$
\widetilde x:=(1_G, \dots,1_G,c_{n+1},c_{n+2},\dots).
$$
Of course, $\widetilde x\in [1_G]_0$ and  $(x,\widetilde x)\in\Cal R$.
In view of the equality in the lefthand side of~\thetag{1-9}, there exist $a_j\in A_j$
and $b_j\in B_j$ such that $c_j=a_jb_j$ for all $j\in\Bbb N$.
Moreover, $a_j=b_j=1_G$ if $j\le n$ according to \thetag{1-11} and~\thetag{1-12}.
Then
$$
\phi(\widetilde x)=(1_G,\dots,1_G, a_{n+1}, b_{n+1}a_{n+2}, b_{n+2}a_{n+3}, \dots)\in[1_G]_0\subset X'.
$$
  It follows from the equality in the righthand side of  \thetag{1-10} that
   $f_na_{n+1}\in F_{n+1}'$.
Hence, the point $y:=(f_na_{n+1},b_{n+1}a_{n+2}, b_{n+2}a_{n+3}, \dots)$ belongs to $X_{n+1}'\subset X'$ and
$(y,\phi(\widetilde x))\in\Cal R'$.

\endcomment %%%%%%%%%%%%%%%%%%%%%%%%%%%%%

We now define a mapping $\psi_{\boldsymbol A,\boldsymbol B}:X\to X'$.
Given $x\in X$, we find $n\ge 0$ such that $x=(f_n,c_{n+1}, c_{n+2}, \dots)\in X_n\subset X$. 
It follows from the lefthand side of \thetag{1-9} and~\thetag{1-11} that there exist unique
$a_{j-1}\in A_{j-1}$ and $b_j\in B_j$ such that $c_{j}=a_{j-1}b_j$ for each $j>n$.
Then we let
$$
\psi_{\boldsymbol A,\boldsymbol B}(x):=(f_na_{n}, b_{n+1}a_{n+1}, b_{n+2}a_{n+2},\dots).
$$
It follows from the righthand sides of \thetag{1-9} and \thetag{1-10} that
$\psi_{\boldsymbol A,\boldsymbol B}(x)\in X_{n}'\subset X'$.
A straightforward verification gives that $\psi_{\boldsymbol A,\boldsymbol B}(x)$
is well defined (i.e. does not depend on the choice of  $n$ such that $x\in X_n$).

Let $\Cal R$ and $\Cal R'$ denote the tail equivalence relations on $X$ and $X'$
respectively.
Let $\alpha$ and $\alpha'$ stand for the return time cocycles of 
$\Cal R$ and $\Cal R'$ respectively.

\proclaim{Proposition 1.6}
\roster
%\item"(i)" $\psi_{\boldsymbol A,\boldsymbol B}$ is a well defined mapping from $X$ to $X'$.
%\item"(ii)" $\psi_{\boldsymbol A,\boldsymbol B}\restriction [1_G]_0=\phi$.
% \item"(i)" $\psi_{\boldsymbol A,\boldsymbol B}$ is a bijection.
%   \item"(iv)" $\psi_{\boldsymbol A,\boldsymbol B}$ is onto.
  \item"(i)"  $\psi_{\boldsymbol A,\boldsymbol B} $ is a homeomorphism of $X$ onto $X'$.
   \item"(ii)"  $\psi_{\boldsymbol A,\boldsymbol B} $ maps the $
   \Cal R$-class of each $x\in X$ bijectively onto the  $\Cal R'$-class of 
   $\psi_{\boldsymbol A,\boldsymbol B}(x) $.
   \item"(iii)" 
    $\psi_{\boldsymbol A,\boldsymbol B} $ transfers the Haar measure on $X$
    to the Haar measure on $X'$.
 \item"(iv)" $\alpha'(\psi_{\boldsymbol A,\boldsymbol B}(x),\psi_{\boldsymbol A,\boldsymbol B}(\widetilde x))=\alpha(x,\widetilde x)$ for all $(x,\widetilde x)\in\Cal R$.
\endroster
\endproclaim

\demo{Proof} 
(i) We first prove that    $\psi_{\boldsymbol A,\boldsymbol B}$ is one-to-one.
If    $\psi_{\boldsymbol A,\boldsymbol B}(x)=\psi_{\boldsymbol A,\boldsymbol B}(\widetilde x)$
for some $x,\widetilde x\in X$, we can find $n\in\Bbb N$ and  elements $f_n,\widetilde f_n\in F_n$, $a_{j-1}\in A_{j-1}$ and $b_j\in B_j$  for all $j>n$ such that
$$
x=(f_n, a_{n}b_{n+1}, a_{n+1}b_{n+2},\dots)\quad\text{and}\quad
\widetilde x=(\widetilde f_n, \widetilde a_{n} \widetilde b_{n+1}, \widetilde a_{n+1}\widetilde b_{n+2},\dots).
$$
Since $\psi_{\boldsymbol A,\boldsymbol B}(x)=\psi_{\boldsymbol A,\boldsymbol B}(\widetilde x)$, it follows that $f_na_{n}=\widetilde f_n \widetilde a_{n}$ and 
$ b_{j}a_{j}=\widetilde b_{j} \widetilde a_{j}$ for each 
$j>n$.
Then~\thetag{1-11} yields that $a_{j}=\widetilde a_{j}$ and $b_{j}=\widetilde b_{j}$ 
for all $j>n$.
We obtain that
$$
f_{n+1}=f_na_{n}b_{n+1}=\widetilde f_n \widetilde a_{n}\widetilde b_{n+1}=\widetilde f_{n+1}.
$$
Hence,
$$
x=(f_{n+1}, a_{n+1}b_{n+2},a_{n+2}b_{n+3},\dots)=(\widetilde f_{n+1}, \widetilde a_{n+1}\widetilde b_{n+2},\widetilde a_{n+2}\widetilde b_{n+3},\dots)=\widetilde x,
$$
as desired.

We now show that  $\psi_{\boldsymbol A,\boldsymbol B}$ is onto.
Take $x'\in X'$ and find $n\in\Bbb N$ such that $x'\in X_n'$.
Then there exist $f_n'\in F_n'$ and $b_j\in B_j$ and  $a_{j+1}\in A_{j+1}$ for each $j>n$ such that
$x'=(f_n', b_{n+1}a_{n+1}, b_{n+2}a_{n+2},\dots)$.
By the lefthand side of \thetag{1-10}, the element $f_{n+1}:=f_n'b_{n+1}$ belongs to
$F_{n+1}$.
It is straightforward to verify that
$$
\psi_{\boldsymbol A,\boldsymbol B}(f_{n+1}, a_{n+1}b_{n+2}, a_{n+2}b_{n+3},\dots)=x'.
$$
 It is easy to see that $\psi_{\boldsymbol A,\boldsymbol B}$ is continuous.
 The  $\psi_{\boldsymbol A,\boldsymbol B}$-image of a cylinder in $X$ is a cylinder in $X'$.
Hence, the mapping  $\psi_{\boldsymbol A,\boldsymbol B}^{-1}$ is also  continuous.
 Thus, (i) is proved.

(ii)--(iv) are routine.
\qed
\enddemo

\definition{Definition 1.7}
We call $\psi_{\boldsymbol A,\boldsymbol B}$ the  {\it  $(\boldsymbol A,\boldsymbol B)$-chain equivalence} of $X$ onto  $X'$.
\enddefinition

Let $\Cal T$ and $\Cal T'$ both satisfy  \thetag{1-4} or \thetag{1-5}. 
Then  the $(C,F)$-actions $T$ and $ T'$ associated with $\Cal T$ and $\Cal T'$ respectively
are well defined.
If $\Cal T$ is chain equivalent to  $\Cal T'$   %and \thetag{1-10} holds
 then
it follows from Proposition~1.6(iv) that the chain equivalence intertwines
$T$ with $T'$, i.e. 
$$
\psi_{\boldsymbol A,\boldsymbol B}\circ T_g=T_g'\circ\psi_{\boldsymbol A,\boldsymbol B}
\quad\text{for all $g\in G$.}
$$
%\comment

It is easy to verify that if %$\Cal T$ and $\Cal T'$
%satisfy \thetag{1-7}, 
 %$\Cal T$ is congruent to $\Cal T'$ and  
$1_G\in \bigcap_{n=0}^\infty (F_n\cap F_n')$ then for each $n\ge 1$,
$$
\aligned
&\psi_{\boldsymbol A,\boldsymbol B}([1_G]_{n-1})=[A_{n-1}]_{n-1} \quad\text{and}\quad
\psi_{\boldsymbol A,\boldsymbol B}^{-1}([1_G]_{n-1})=[B_{n}]_{n}. % \\
%&\psi_{\boldsymbol A,\boldsymbol B}([1_G]_1)=[A_2]_2, \quad
%\psi_{\boldsymbol A,\boldsymbol B}^{-1}([1_G]_2)=[B_2]_2,\\
%&\psi_{\boldsymbol A,\boldsymbol B}([1_G]_2)=[A_3]_3, \quad
%\psi_{\boldsymbol A,\boldsymbol B}^{-1}([1_G]_3)=[B_3]_3,\\
%&\dots
\endaligned
\tag1-12
$$
Thus, \thetag{1-12} gives formulae for how to ``reconstruct'' the sequences
$\boldsymbol A$ and $\boldsymbol B$ if $\psi_{\boldsymbol A,\boldsymbol B}$
is known.

\remark{Remark 1.8}
Suppose now that  $\Cal T$ is chain equivalent  to $\Cal T'$ which satisfies \thetag{1-7}.
Then there is a calibration $\widetilde{\Cal T}$ of $\Cal T$ such that:
\roster
\item"(i)"
$\widetilde {\Cal T}$ satisfies \thetag{1-7},
\item"(ii)"
$\widetilde{\Cal T}$  is chain equivalent  to  $\Cal T'$ and
if  $(\widetilde{A_n})_{n=1}^\infty$ and $(\widetilde{B_n})_{n=1}^\infty$ stand for  the corresponding sequences
of finite subsets in $G$ (satisfying \thetag{1-9}--\thetag{1-11}) then
$1\in \bigcap_{n=1}^\infty (\widetilde A_n\cap\widetilde B_n)$.
\endroster
Indeed,
it follows from the right equation in~\thetag{1-9} and \thetag{1-11} that for each $n\in\Bbb N$, there exist unique $b_n\in B_n$ and $a_{n}\in A_{n}$
such that $b_na_{n}=1_G$.
Hence, $a_{n}^{-1}=b_n$.
Let $\widetilde B_n:=B_nb_n^{-1}$, $\widetilde A_{n}:=b_nA_{n}$ and  $\widetilde F_{n}:=F_{n}b_n^{-1}$.
We also let $b_0:=\{1_G\}$ and $\widetilde A_0:=b_0 A_0=A_0.$
Since $F_0'=\{1_G\}$, it follows from \thetag{1-10} that $F_0A_0=\{1_G\}$.
We now deduce  from \thetag{1-9}--\thetag{1-11} that the following hold for each $n\in\Bbb N$:
$$
\gathered
\widetilde A_{n-1}\widetilde B_n=b_{n-1}C_nb_n^{-1},\quad \widetilde B_n \widetilde A_{n}=B_nA_n=C'_{n},\\
F'_{n-1}\widetilde B_nb_n\subset F_n,\quad  F_{n-1}b_{n-1}^{-1}\widetilde A_{n-1}\subset F'_{n-1}\quad\text{and}
\\
\widetilde A_{n-1}^{-1}\widetilde A_{n-1}\cap \widetilde B_n\widetilde B_n^{-1}=\widetilde B_n^{-1}\widetilde B_n\cap
\widetilde A_{n}\widetilde A_{n}^{-1}=\{1_G\}.
\endgathered
\tag1-13
$$
Let $\boldsymbol z:=(b_0^{-1},b_1^{-1},b_2^{-1},\dots)$.
Denote by $\widetilde{\Cal T}$ the $\boldsymbol z$-calibration of $\Cal T$.
Of course,  $\widetilde{\Cal T}$ satisfies \thetag{1-7} and $1\in \bigcap_{n=1}^\infty (\widetilde A_n\cap\widetilde B_n)$.
It follows from \thetag{1-13} that   $\widetilde{\Cal T}$ is chain equivalent  to ${\Cal T'}$.

\endremark

%\endcomment

\proclaim{Fact E \cite{Da3}}
Let $\Cal T$ and $\Cal T'$ satisfy~\thetag{1-1}, %\thetag{1-3} 
and \thetag{1-4}.
Then the topological $(C,F)$-actions of $G$ associated with $\Cal T$ and $\Cal T'$ are topologically isomorphic if and only if there exist two sequences
$\boldsymbol k:=(k_n)_{n=1}^\infty$ and 
$\boldsymbol l:=(l_n)_{n=1}^\infty$ of nonnegative integers such that
$0=k_0=l_0<k_1<l_1<k_2<l_2<\cdots$  and
the $\boldsymbol k$-telescoping of $\Cal T$ is chain equivalent  to the
$\boldsymbol l$-telescoping of $\Cal T'$.
\endproclaim

We will also utilize  the following fact to prove the main results of the paper.

\proclaim{Proposition 1.9}
Let $T$ be the $(C,F)$-action associated with a $(C,F)$-sequence $\Cal T$ satisfying
\thetag{1-1} %\thetag{1-3} 
 and \thetag{1-5}.
Denote by $(X,\mu)$ the corresponding $(C,F)$-space endowed with the Haar measure.
Let $Q=(Q_g)_{g\in G}$ be an ergodic measure preserving action of $G$ on 
a $\sigma$-finite standard measure space $(Y,\nu)$.
% and $\mu(X)=\nu(Y)$.
Let $\phi,\psi:(X,\mu)\to (Y,\nu)$  be two measure preserving isomorphisms
such that $\phi T_g\phi^{-1}=\psi T_g\psi^{-1}=Q_g$ for each $g\in G$.
If 
$$
\lim_{n\to\infty}\frac{\nu(\phi([1_G]_n)\triangle\psi([1_G]_n))}{\mu([1_G]_n)}= 0
$$
then $\phi=\psi$ almost everywhere.
\endproclaim

\demo{Proof}
 We note that for each subset $A\subset F_n$,
$$
\align
\nu(\phi([A]_n)\triangle\psi([A]_n))&=
\nu\Bigg(\bigg(\bigsqcup_{g\in A}\phi(T_g[1_G]_n)\bigg)\triangle\bigg(\bigsqcup_{g\in A}\psi(T_g[1_G]_n)\bigg)\Bigg)\\
&\le\sum_{g\in A}\nu\big(Q_g\phi([1_G]_n)\triangle Q_g \psi([1_G]_n)\big)\\
&\le \#A\cdot\nu(\phi([1_G]_n)\triangle\psi([1_G]_n)).
\endalign
$$
Hence, 
$$
\lim_{n\to\infty}\max_{A\subset F_n}
\frac{\nu(\phi([A]_n)\triangle\psi([A]_n))}{\mu([A]_n)}=
\lim_{n\to\infty}\max_{A\subset F_n}
\frac{\nu(\phi([A]_n)\triangle\psi([A]_n))}{\# A\cdot \mu([1_G]_n)}= 0.
$$
It follows that $\nu(\phi(B)\triangle\psi(B))=0$ for each cylinder $B$ in $X$.
This implies that 
$\nu(\phi(D)\triangle\psi(D))=0$ for each
 Borel subset $D\subset X$ of finite measure.
Hence, $\phi=\psi$ almost everywhere.
\qed
\enddemo

\comment

From Proposition~1.9 and \thetag{1-13} we deduce the following.

\proclaim{Corollary 1.10} Let $\Cal T$ and $\Cal T'$ be congruent
$(C,F)$-sequences and let  \thetag{1-5} and \thetag{1-12} hold.
Denote by $T$ and $T'$  the $(C,F)$-actions associated with $\Cal T$ and $\Cal T'$ respectively.
Let $(X,\mu)$ and $(X',\mu')$ denote the corresponding $(C,F)$-spaces. 
Let $\phi:X\to X'$ be a measure preserving isomorphism of $T$ onto $T'$.
If 
$$
\lim_{n\to\infty}\frac{\mu'(\phi([1_G]_n)\triangle [A_{n+1}]_{n+1})}{\mu([1_G]_n)}= 0
$$
then $\phi=\psi_{\boldsymbol A,\boldsymbol B}$ a.e.
\endproclaim

\endcomment

\comment

Let $S=(S_g)_{g\in G}$
and $R=(R_g)_{g\in G}$ be ergodic measure preserving free actions on standard finite measure spaces
$(Y,\nu)$ and $(Z,\xi)$ respectively.
Let $A\subset Y$ and $B\subset Z$ be subsets of positive measure.
Suppose that there is a Borel isomorphism $\psi:A\to B$ such that if $a\in A$ then
 $S_ga\in A$ for some $g\in  G$ if and only if  $\psi(S_ga)\in B$ and $\psi(S_ga)=R_g\psi(a)$.

Let $\Cal T$ and $\Cal T'$ be two $(C,F)$-sequences.
If
there is 
an
 increasing sequence of integers
 $$0=l_0<l_1'<l_1<l_2'<l_2<\cdots$$ 
 and subsets
$A_n \subset  F_{l_n'}'$,  $B_n \subset  F_{l_n}$, 
such that
$$
\gathered
A_{n}B_{n}=C_{l_{n-1}+1}\cdots C_{l_{n}},\quad
B_{n} A_{n+1}=C'_{l_{n}'+1}\cdots C'_{l_{n+1}'},
\\
F'_{l_n'}B_n\subset F_{l_n},\quad F_{l_n}A_{n+1}\subset F'_{l_{n+1}},
\\
(F'_{l_n'})^{-1}F'_{l_n'}\cap B_nB_n^{-1}=F_{l_n}^{-1}F_{l_n}\cap A_{n+1}A_{n+1}^{-1}=\{1\}
\endgathered
\tag1-7
$$
 for each $n>0$.

$(\Leftarrow)$
Suppose that  \thetag{2-1} is satisfied.
Let $\widetilde T$ and $\widetilde T'$ denote the $(C,F)$-actions associated with  the $(l_{i})_{i\ge 0}$-telescoping of $(C_n,F_{n-1})_{n\ge 1}$ and the $(l'_{i})_{i\ge 1}$-telescoping of 
  $(C_n',F_{n-1}')_{n\ge 1}$ respectively.
  Let $\widetilde X$ and $\widetilde X'$ denote the $(C,F)$-spaces of these actions.
  Of course, $\widetilde T$ is isomorphic to $T$ and 
  $\widetilde T'$ is isomorphic to $T'$.
  Take $x\in [1]_0\subset\widetilde X$.
In view of \thetag{2-1}, we have a unique  expansion  $x=(1,a_1b_1,a_2b_2,\dots)\in[1_G]_0$ of 
$x$ such that $a_i\in A_i$ and $b_i\in B_i$ for all $i\ge 1$.
We now set
$$
\phi(x):=(a_1, b_1a_2,b_2a_3,\dots)\in[A_1]_1\subset \widetilde X'.
$$
It is standard to verify that  $\phi$ is a homeomorphism of  the cylinder $[1]_0\subset \widetilde X$ onto the cylinder $[A_1]_1\subset\widetilde X'$.
 Then 
 $$
 \align
 \{g\in G\mid  \widetilde T_gx\in [1]_0\}&=\bigcup_{n=1}^\infty  (A_1 B_1)\cdots (A_n B_n)(a_nb_n)^{-1}\cdots(a_1b_1)^{-1}\\
&=\bigcup_{n=1}^\infty A_1 (B_1A_2)\cdots (B_n A_{n+1})(b_na_{n+1})^{-1}\cdots(b_1a_2)^{-1}a_1^{-1}\\
& =\{g\in G\mid  \widetilde T_g'\phi(x)\in [A_1]_1\}.
 \endalign
 $$
The equality $ \phi(\widetilde T_gx)=\widetilde T_g\phi(x)$ whenever $\widetilde T_gx\in[1]_0$ is now verified in a straightforward way.
 It remains to apply Lemma~2.1.

We call this isomorphism {\it the braid isomorphism}.

\endcomment

\head 2. Isomorphic $(C,F)$-actions
\endhead

Let  $\Cal T=(C_n, F_{n-1})_{n>0}$ and $\widetilde{\Cal  T}=(\widetilde C_n, \widetilde F_{n-1})_{n>0}$
be two $(C,F)$-sequences satisfying 
\thetag{1-1}, %\thetag{1-3}, 
\thetag{1-5} and \thetag{1-7}.
Denote by $T=(T_g)_{g\in G}$ and $\widetilde T=(\widetilde T_g)_{g\in G}$ the $(C,F)$-actions associated with $\Cal T$ and $\widetilde{\Cal T}$ respectively.
Then $T$ and $\widetilde T$ are well defined  measure preserving actions on standard non-atomic $\sigma$-finite spaces
$(X,\goth B,\mu)$ and $(\widetilde X,\widetilde{\goth B},\widetilde\mu)$ respectively.
Here $X$ and $\widetilde X$ are the $(C,F)$-spaces associated with $\Cal T$ and $\Cal T'$ respectively.
If $\Cal T$ and $\widetilde{\Cal T}$ satisfy \thetag{1-3}
then $\mu$ and $\widetilde\mu$ will denote 
 the normalized (i.e. probability) Haar measures for the tail equivalence relations on $X$ and $\widetilde X$ respectively.

We will use the following notation below: if $n>m$ then we denote by $C_{n,m}$
the product $C_nC_{n+1}\cdots C_m\subset G$.
The product  $\widetilde C_{n,m}$ is defined in a similar way.
%We will assume that $\mu(X)=\widetilde\mu(\widetilde X)=1$.
Fix a decreasing sequence $(\epsilon_n)_{n=1}^\infty$ of positive reals such that  $\epsilon_1>1$ and $\sum_{n=2}^\infty\epsilon_n<\frac13$.

We now state and prove  the first main result of the paper.

\proclaim{Theorem 2.1} The $(C,F)$-actions $T$ and $\widetilde T$ are measure theoretically isomorphic if and only if there exist a sequence
$$
0=k_0=l_0=k_1<l_1<k_2<l_2<\cdots
$$
 of non-negative integers and subsets $J_n\subset F_{k_{n}}$, $\widetilde J_n\subset \widetilde F_{l_{n}}$ such that
\roster
\item"(i)"
$F_{k_{n}} \widetilde J_n \subset \widetilde F_{l_{n}}$,
\item"(ii)" 
the mapping $F_{k_{n}}\times \widetilde J_n \ni (f, \widetilde f) \mapsto f\widetilde f \in \widetilde F_{l_{n}}$ is one-to-one, 
\item"(iii)" 
%$\dfrac{\# \widetilde F_{l_{n}} - \# F_{k_{n}} \# \widetilde J_n}{\# \widetilde F_{l_{n}}} < 2\epsilon_n$,
%\item"(iv)" 
$\dfrac{\#((\widetilde J_n J_{n+1}) \Delta C_{k_{n} + 1, k_{n+1}})}{\# C_{k_{n} + 1, k_{n+1}}} < 2\epsilon_n$,
\item"(i)$'$" $\widetilde F_{l_{n}}  J_{n+1} \subset F_{k_{n+1}}$,
\item"(ii)$'$" 
the mapping $\widetilde F_{l_{n}}\times J_{n+1} \ni (\widetilde f, f) \mapsto \widetilde f f \in  F_{k_{n+1}}$ is one-to-one,
%(i.e. $(\widetilde F_{l_{n+1}} - \widetilde F_{l_{n+1}}) \cap (J_{n+1} - J_{n+1}) = \{0\}$),
\item"(iii)$'$"
% $\dfrac{\# F_{k_{n+1}} - \# \widetilde F_{l_{n}} \# J_{n+1}}{\# F_{k_{n+1}}} < 2\epsilon_n$,
%\item"(iv)$'$" 
$\dfrac{\#((J_{n+1}  \widetilde J_{n+1}) \Delta \widetilde C_{l_n + 1, l_{n + 1}})}{\# \widetilde C_{l_n + 1, l_{n+1}}} < 2\epsilon_n$
\endroster
for each $n\ge 0$.
\endproclaim

\demo{Proof} We first prove the ``only if" claim.
 Let $\phi:(X,\goth B,\mu)\to (\widetilde X,\widetilde{\goth B},\widetilde\mu)$ be a measure preserving isomorphism\footnote{If $T$ and $\widetilde T$ are isomorphic via a nonsingular isomorphism $\phi$, i.e. $\phi T_g\phi^{-1}=\widetilde T_g$ for each $g\in G$ and $\mu\sim\widetilde\mu\circ\phi$ then  it is easy to verify that $\mu=d\cdot\widetilde \mu\circ\phi$ for some constant $d>0$. 
 In this case we replace $\widetilde\mu$ with $d\cdot\widetilde \mu$.
 Then $\phi$ will be measure preserving.} 
 that intertwines
$T$ with $\widetilde T$.
We will construct the desired objects via an inductive process.
On the first step we  let $\widetilde J_0:=\{1_G\}$ and $J_{1}:=\{1_G\}$.
Suppose that for some $n\in\Bbb N$, we have already constructed  a finite sequence of integers $0<l_1<k_2<l_2<\cdots<k_n$  and 
subsets $(J_m)_{m=1}^n$ and $(\widetilde J_m)_{m=1}^{n-1}$  that satisfy (i)--(iii) 
and (i)$'$--(iii)$'$.
\comment

We will assume, in addition, that the following inequalities hold:
\roster
\item"(v)" $\mu(X\setminus X_{k_n}) < \epsilon_{n}$,%\mu(X)$,
\item"(v)$'$"
 $\widetilde\mu(\widetilde X\setminus \widetilde X_{l_{n-1}}) < \epsilon_{n-1}$.%\widetilde\mu(\widetilde X)$.
\endroster

\endcomment
Our purpose is to find integers $l_n$ and $k_{n+1}$  such that $k_{n+1}>l_n>k_n$ and subsets $\widetilde J_{n}\subset\widetilde F_{l_n}$
and $J_{{n+1}}\subset F_{k_{n+1}}$  for which (i)--(iii) and (i)$'$--(iii)$'$ are satisfied.
Given  $l>k_n$,
we let 
$$
 \widetilde F_{l}^\circ:=\{f\in \widetilde F_{l}\mid F_{k_n}f\subset \widetilde F_{l}\}.
 $$
% According to Fact~B,  $(\widetilde F_l)_{l=1}^\infty$ is a F{\o}lner sequence.
% Hence,  $\# \widetilde F_{l}^\circ/\#\widetilde F_{l}\to 1$ as $l\to\infty$.
The sequence of  rings (of cylinders) $\{[\widetilde A]_{l}\mid \widetilde A\subset\widetilde F_l\}$   approximates the entire Borel $\sigma$-algebra $\widetilde{\goth B}$ as $l\to\infty$.
We claim that the sequence  $\{[\widetilde A]_{l}\mid \widetilde A\subset\widetilde F_l^\circ\}$   also approximates $\widetilde{\goth B}$ (mod $\widetilde\mu$) as $l\to\infty$.
Indeed, take a cylinder $D$ in $\widetilde X$.
Then $D=[\widetilde D]_m$ for some $m\in\Bbb N$ and a subset $\widetilde D\subset \widetilde F_m$.
For each  $l>m$, we have that $D=[\widetilde D\widetilde C_{m+1}\cdots \widetilde C_l]_l$.
For a fixed $n$, we let
$$
\widetilde D^\circ_l:=\{f\in \widetilde D\widetilde C_{m+1}\cdots \widetilde C_l\mid F_{k_n}f\subset \widetilde F_l \}.
$$ 
Of course, $[\widetilde D^\circ_l]_l\subset D$.
It follows from \thetag{1-5} that $\widetilde\mu(D\setminus[\widetilde D^\circ_l]_l)\to 0$ as $l\to\infty$.
Since $D$ is an arbitrary cylinder in $\widetilde X$, it follows that
$\{[\widetilde A]_{l}\mid \widetilde A\subset\widetilde F_l^\circ\}$   approximates the entire Borel $\sigma$-algebra $\widetilde{\goth B}$ as $l\to\infty$, as claimed.

Since $\phi$ is an isomorphism, it
  follows that the sequence of  rings $\{\phi^{-1}([\widetilde A]_{l})\mid \widetilde A\subset\widetilde F_l^\circ\}$   approximates the Borel $\sigma$-algebra $\goth B$ on $X$ (mod $\mu$) as $l\to\infty$.
Hence,
 we
can find $l_n>k_n$ 
  and a subset $\widetilde J_n\subset \widetilde F_{l_n}^\circ$ such that 
$$
\align
%\widetilde\mu(\widetilde X\setminus \widetilde X_{l_n}) &< \epsilon_n\quad\text{and}\tag2-1\\
\mu\big([1_G]_{k_n}\triangle \phi^{-1}( [\widetilde J_n]_{l_n})\big)
&<\epsilon_n\mu([1_G]_{k_n})\tag2-1
\endalign
$$
Since $\widetilde J_n\subset \widetilde F_{l_n}^\circ$, it follows that (i) holds.
%Of course, \thetag{2-1} is exactly (v)$'$ with $n$ in place of $n-1$.
 Since $[\widetilde J_n]_{l_n}=\bigsqcup_{\widetilde f\in\widetilde J_n}[\widetilde f]_{l_n}$ and \thetag{2-1} holds, we can  assume without loss of generality (passing, if necessarily, to a subset in $\widetilde J_n$)   that
 $$
\mu\big([1_G]_{k_n}\cap \phi^{-1}( [\widetilde f]_{l_n})\big)>0.5\mu(\phi^{-1}( [\widetilde f]_{l_n}))\quad\text{for each $\widetilde f\in \widetilde J_n$.
}\tag2-2
$$
For each $f\in F_{k_n}$, we have that $f\widetilde J_n\subset \widetilde F_{l_n}$ and hence
$\widetilde T_f[\widetilde J_n]_{l_n}=[f\widetilde J_n]_{l_n}$.
Therefore,
$$
\mu\big([1_G]_{k_n}\cap \phi^{-1}( [\widetilde J_n]_{l_n})\big)=
\mu\big(T_f[1_G]_{k_n}\cap \phi^{-1}(\widetilde T_f[\widetilde J_n]_{l_n})\big)=
\mu\big([f]_{k_n}\cap \phi^{-1}([f\widetilde J_n]_{l_n})\big).
$$
Hence, we deduce from \thetag{2-1} and \thetag{2-2} that 
$$
\align
\mu\big([f]_{k_n}\cap \phi^{-1}( [f\widetilde J_n]_{l_n})\big)&>(1-\epsilon_n)\mu([f]_{k_n})\quad\text{and}\tag2-3\\
\mu\big([f]_{k_n}\cap \phi^{-1}( [f\widetilde f]_{l_n})\big)&>0.5\mu(\phi^{-1}( [f\widetilde f]_{l_n}))\tag2-4
\endalign
$$
for each $\widetilde f\in \widetilde J_n$.
Since the cylinders $[f]_{k_n}$, $f\in F_{k_n}$, are mutually disjoint, it follows from \thetag{2-4} that the subsets $f\widetilde J_n$, $f\in F_{k_n}$, are mutually disjoint.
Thus, (ii) holds. 
\comment

We note that~\thetag{2-4} and the equality $X_{k_n}=\bigsqcup_{f\in F_{k_n}}[f]_{k_n}$ yield
the following:
$$
\align
\widetilde\mu\bigg(\phi(X_{k_n})\cap\bigg(\bigsqcup_{f\in F_{k_n}}[f\widetilde J_n]_{l_n})\bigg)&=
\mu\Bigg(X_{k_n}\cap\bigg(\bigsqcup_{f\in F_{k_n}}\phi^{-1}([f\widetilde J_n]_{l_n})\bigg)\Bigg)\\
&\ge\sum_{f\in F_{k_n}}\mu([f]_{k_n}\cap \phi^{-1}( [f\widetilde J_n]_{l_n}))\\
&>(1-\epsilon_n)\mu(X_{k_n}).
\endalign
$$
Using this inequality and (v), we obtain that
$$
\align
\widetilde\mu\big([F_{k_n}\widetilde J_n]_{l_n}\big)&\ge \widetilde\mu(\phi( X_{k_n})\cap[F_{k_n}\widetilde J_n]_{l_n})\\
&>(1-\epsilon_n)\mu(X_{k_n})\\
&>(1-\epsilon_n)^2\\
&>(1-\epsilon_n)^2\widetilde\mu(\widetilde X_{l_n}).
\endalign
$$
Hence $\frac{\# F_{k_n}\# \widetilde J_n}{\#\widetilde F_{l_n}}>(1-\epsilon_n)^2>1-2\epsilon_n$ and (iii) follows.

\endcomment

Arguing in a similar way, we can find $k_{n+1}>l_n$ and a subset 
$$
J_{n+1}\subset \{f\in F_{k_{n+1}}\mid \widetilde F_{l_n}f\subset F_{k_{n+1}}\}
$$
 such that
$$
\gather
%\mu( X\setminus  X_{k_{n+1}}) < \epsilon_{n+1}\quad\text{and}\\
%\# \widetilde F_{l_n}^\circ>(1-\epsilon_n/\# F_{k_n}) \# \widetilde F_{l_n}\quad\text{and}\tag2-2\\
\widetilde\mu\big([1_G]_{l_n}\triangle \phi( [J_{n+1}]_{k_{n+1}})\big)<
%\frac{\epsilon_{n+1}\mu([1_G]_{k_n})}{\# \widetilde F_{l_n}}
\epsilon_{n+1}
\widetilde\mu([1_G]_{l_n}).\tag2-5
\endgather
$$
In turn, this inequality imply (i)$'$--(iii)$'$ in a similar way as \thetag{2-1} 
%and  \thetag{2-2} 
implied~(i)--(iii).
Since $[1_G]_{k_n}=[C_{k_n+1,k_{n+1}}]_{k_{n+1}}$, it follows  from \thetag{2-1} and \thetag{2-5} that
$$
\align
\epsilon_n\mu([1_G]_{k_n}) &>
\mu([1_G]_{k_n}\triangle \phi^{-1}([\widetilde J_n]_{l_n}))\\
&=\mu\bigg([C_{k_n+1,k_{n+1}}]_{k_{n+1}}\triangle \bigsqcup _{g\in \widetilde J_n}T_g\phi^{-1}( [ 1_G]_{l_n})\bigg)
\\
&=
\widetilde\mu\bigg(\phi\big([C_{k_n+1,k_{n+1}}]_{k_{n+1}}\big)\triangle \bigsqcup _{g\in \widetilde J_n}\widetilde T_g [ 1_G]_{l_n}\bigg)
\\
&>
\widetilde\mu\bigg(\phi\big([C_{k_n+1,k_{n+1}}]_{k_{n+1}}\big)\triangle \bigsqcup _{g\in \widetilde J_n}
\widetilde  T_g\phi([J_{n+1}]_{k_{n+1}})\bigg)\\
%\frac{\#\widetilde J_n\epsilon_{n+1}\mu([1_G]_{k_n})\widetilde\mu([1_G]_{l_n})}
%{\# \widetilde F_{l_n}}\\
&-\#\widetilde J_n\widetilde\mu\big([1_G]_{l_n}\triangle\phi([J_{n+1}]_{k_{n+1}})\big)\\
&>
\mu\big([C_{k_n+1,k_{n+1}}]_{k_{n+1}}\triangle [\widetilde J_n J_{n+1}]_{k_{n+1}})
-\frac{\epsilon_{n+1}}{1-\epsilon_{n+1}}\#\widetilde J_n\widetilde\mu(\phi([J_{n+1}]_{k_{n+1}}))\\
&=
\mu\big([C_{k_n+1,k_{n+1}}\triangle (\widetilde J_n J_{n+1})]_{k_{n+1}})-\frac32\epsilon_{n+1}\# \widetilde J_n\# J_{n+1}\mu([1_G]_{k_{n+1}})\\
&=
\Big(\#(C_{k_n+1,k_{n+1}}\triangle (\widetilde J_n J_{n+1}))-
\frac32\epsilon_{n+1}\#(\widetilde J_nJ_{n+1})\Big)
\mu([1_G]_{k_{n+1}})
.
\endalign
$$
Hence
$$
\epsilon_n >\frac{\#(C_{k_n+1,k_{n+1}}\triangle (\widetilde J_n J_{n+1}))}
{\# C_{k_n+1,k_{n+1}}}
-\frac32\epsilon_{n+1}\frac{\#(\widetilde J_n J_{n+1})}{\# C_{k_n+1,k_{n+1}}}.
$$
This yields that
$$
\frac{\#(C_{k_n+1,k_{n+1}}\triangle (\widetilde J_n J_{n+1}))}
{\# C_{k_n+1,k_{n+1}}}
<
\frac{\epsilon_n+1.5\epsilon_{n+1}}{1-1.5\epsilon_{n+1}}
$$
and (iii) follows.
The inequality (iii)$'$ is proved is a similar way.
We start with the inequality which is \thetag{2-5} but with $n-1$ in place of $n$.
Without loss of generality, we may assume that this inequality holds by the inductive assumption.
Then we have:
$$
\align
\epsilon_{n-1}&\widetilde\mu([1_G]_{l_{n-1}})>
\widetilde\mu([1_G]_{l_{n-1}}\triangle \phi([J_n]_{k_n}))\\
&=
\mu\bigg(\phi^{-1}\big([\widetilde C_{l_{n-1}+1,l_n}]_{l_{n}}\big)\triangle \bigsqcup _{g\in J_n}T_g [ 1_G]_{k_n}\bigg)
\\
&>
\mu\bigg(\phi^{-1}\big([\widetilde C_{l_{n-1}+1,l_n}]_{l_{n}}\big)\triangle \bigsqcup _{g\in J_n}T_g \phi^{-1}([ \widetilde J_n]_{l_n})\bigg)
-\# J_n\mu\big([1_G]_{k_n}\triangle\phi^{-1}[\widetilde J_n]_{l_n}\big)\\
&>
\widetilde\mu\big([\widetilde C_{l_{n-1}+1,l_n}]_{l_{n}}\triangle [J_n\widetilde J_n]_{l_n}\big)
-
\frac{\epsilon_n}{1-\epsilon_n}\# J_n\mu(\phi^{-1}[\widetilde J_n]_{l_n}).
\endalign
$$
Hence,
$$
\epsilon_{n-1}>\frac{\#(\widetilde C_{l_{n-1}+1,l_n}\triangle(J_n\widetilde J_n))}{\# \widetilde C_{l_{n-1}+1,l_n}}
-\frac{3\epsilon_n\# (J_n\widetilde J_n)}{2\#\widetilde C_{l_{n-1}+1,l_n}},
$$
which implies (iii)$'$.
% on the next inductive step.
Thus, the ``only if'' part of the theorem is proved.

We now prove the ``if'' claim.
Thus, suppose that 
there exist  a sequence 
$$
0=k_0=l_0=k_1<l_1<k_2<l_2<\cdots
$$
 of integers and subsets $J_n\subset F_{k_{n}}$, $\widetilde J_n\subset \widetilde F_{l_{n}}$ such that
(i)--(iii) and (i)$'$--(iii)$'$ are satisfied.
Consider two  sequences
$$
\Cal V:=(F_{k_n}, \widetilde J_nJ_{n+1})_{n=0}^\infty\quad\text{and}\quad
\Cal W:=(\widetilde F_{l_n}, J_{n+1}\widetilde J_{n+1})_{n=0}^\infty
$$ 
 of finite subsets in $G$.
Of course, $\#(\widetilde J_nJ_{n+1})>1$.
It follows from (i) and (i)$'$ that  $F_{k_n} \widetilde J_nJ_{n+1}\subset F_{k_{n+1}}$.
We deduce from  (ii) and (ii)$'$ that  $F_{k_n}c\cap F_{k_n} c'=\emptyset$ for all $c,c'\in \widetilde J_nJ_{n+1}$ if $c\ne c'$.
Thus, $\Cal V$ is a $(C,F)$-sequence that satisfies \thetag{1-1}.
In a similar way,   one can verify that $\Cal W$  is also a $(C,F)$-sequence that satisfies \thetag{1-1}.

We claim  that $\Cal V$ is chain equivalent  to  $\Cal W$.
Let $\widetilde{\boldsymbol J}:=(\widetilde J_n)_{n=0}^\infty$ and 
$\boldsymbol J:=( J_{n+1})_{n=1}^\infty$. 
Then~\thetag{1-9} holds for $\Cal V$ and $\Cal W$ by the definition of $\Cal V$ and $\Cal W$ with $\widetilde{\boldsymbol J}$ and  $\boldsymbol J$ 
in place of $\boldsymbol A$ and $\boldsymbol B$ respectively.
The inclusions \thetag{1-10} follow from  (i) and (i)$'$.
Finally,~\thetag{1-11} follow from (ii) and (ii)$'$.
Thus, the claim  is proved.
%?????????????????It follows that the $(C,F)$-actions associated with $\Cal V$ and $\Cal W$ are isomorphic.

We now let 
 $$
 \gather
 A_n:=(\widetilde J_n J_{n+1}) \cap C_{k_n + 1, k_{n+1}},\quad 
 B_n:=( J_{n+1} \widetilde J_{n+1}) \cap \widetilde C_{l_n + 1, l_{n+1}},\\
  \boldsymbol A:=(A_n)_{n=1}^\infty,\quad { \boldsymbol B}:=( B_n)_{n=1}^\infty, \\
 \boldsymbol k:=(k_n)_{n=0}^\infty \quad
 \text{and}\quad  \boldsymbol l:=(l_n)_{n=0}^\infty.
  \endgather
 $$
It follows from (iii) and (iii)$'$ that 
$$
\align
\sum_{n=1}^\infty\bigg(1-\frac{\#A_n}{\# C_{k_n + 1, k_{n+1}}}\bigg) &<2\sum_{n=1}^\infty\epsilon_n<\infty 
\quad \text{and}\tag2-6
\\
\sum_{n=1}^\infty\bigg(1-\frac{\#B_n}{\# \widetilde C_{l_n + 1, l_{n+1}}}\bigg) &<2\sum_{n=1}^\infty\epsilon_n<\infty. \tag2-7
\endalign
$$
In a similar way, 
$$
\align
\sum_{n=1}^\infty\bigg(1-\frac{\# A_n}{\# (\widetilde J_n J_{n+1})}\bigg)&<4\sum_{n=1}^\infty\epsilon_n<\infty
\quad
\text{and} \tag2-8\\
\sum_{n=1}^\infty\bigg(1-\frac{\#B_n}{\#( J_{n+1} \widetilde J_{n+1})}\bigg)&<4\sum_{n=1}^\infty\epsilon_n<\infty.
\tag2-9
\endalign
$$ 
Let  $\Cal T_{\boldsymbol A,\boldsymbol k}$ denote
the $\boldsymbol A$-reduction of the $\boldsymbol k$-telescoping of
$\Cal T$.
It is well defined in view of \thetag{2-6}.
Also, let $\Cal V_{\boldsymbol A}$ stand for 
the $\boldsymbol A$-reduction of %the $\boldsymbol k$-telescoping of
 $\Cal V$.
 The reduction is well defined in view of \thetag{2-8}.
Of course,  $\Cal T_{\boldsymbol A,\boldsymbol k}=\Cal V_{\boldsymbol A}$.
It follows from \S1.5 and \$1.6 that the $(C, F)$-action associated with 
$\Cal T_{\boldsymbol A,\boldsymbol k}$ is well defined and  isomorphic to $T$.
Hence, the $(C, F)$-action associated with 
$\Cal V_{\boldsymbol A}$ is also well defined and  isomorphic to $T$.
This implies, in turn, that the $(C, F)$-action associated with 
$\Cal V$ is well defined and  isomorphic to $T$.

In a similar way,  let $\widetilde {\Cal T}_{{\boldsymbol B},\boldsymbol l}$ denote
the ${\boldsymbol B}$-reduction of the $\boldsymbol l$-telescoping of 
$\widetilde{\Cal T}$.
It is well defined in view of \thetag{2-7}.
Let $ {\Cal W}_{{\boldsymbol B}}$ stand for the ${\boldsymbol B}$-reduction of
$\Cal W$.
It is well defined in view of \thetag{2-9}.
Of course, 
 $\widetilde {\Cal T}_{{\boldsymbol B},\boldsymbol l}= {\Cal W}_{{\boldsymbol B}}$.
 Arguing in the same way as above, we conclude that
 the $(C, F)$-action associated with 
$ {\Cal W}$ is well defined and  isomorphic to $\widetilde T$.

Since $\Cal V$ is chain equivalent  to $\Cal W$, it follows from \S1.7 that the $(C,F)$-actions associated to 
$\Cal V$ and $\Cal W$ are isomorphic. 
Thus, $T$ and $\widetilde T$ are isomorphic, as desired.
 \qed
 \enddemo
 
 We can provide 
 an explicit formula for the  isomorphism $\phi$ between $T$ and $T'$.
 This isomorphism is a composition of 7 mappings, each of which is either  a telescoping (or  inverse to a telescoping),  a reduction (or  inverse to a reduction), or a chain equivalence. 
 We will use below the notation introduced in the statement and the proof of Theorem~2.1.

 \comment
 The $(C,F)$-action associated with $\Cal T_{\boldsymbol A,\boldsymbol k}$ is
 $\rho_{\boldsymbol A}\iota_{\boldsymbol k}T\iota_{\boldsymbol k}^{-1}\rho_{\boldsymbol A}^{-1}$,
 where $\rho_{\boldsymbol A}$ is the $\boldsymbol A$-reduction mapping and
$ \iota_{\boldsymbol k}$
 is the $\boldsymbol k$-telescoping mapping.
 The $(C,F)$-action associated with $\widetilde {\Cal T}_{{\boldsymbol B},\boldsymbol l}$
  is
 $\rho_{\boldsymbol B}\iota_{\boldsymbol l}\widetilde T\iota_{\boldsymbol l}^{-1}\rho_{\boldsymbol B}^{-1}$,
 where $\rho_{\boldsymbol B}$ is the $\boldsymbol B$-reduction mapping and
$ \iota_{\boldsymbol l}$
 is the $\boldsymbol l$-telescoping mapping.
\endcomment

 \comment

 First, since $\Cal V$ and $\Cal W$ are congruent,
 a normalization
 $(\boldsymbol z,\boldsymbol v,\widetilde{\boldsymbol J}',{\boldsymbol J'})$
 for  $(\widetilde{\boldsymbol J},\boldsymbol J)$ is well defined by Lemma~1.6.
 Denote by $\psi_{\widetilde{\boldsymbol J}',{\boldsymbol J}}$ the $(\widetilde{\boldsymbol J}',{\boldsymbol J'})$-congruency mapping.
 Let $\phi_{\boldsymbol z}$ and $\phi_{\boldsymbol v}$ stand for the
 corresponding callibration mappings.

 \endcomment
 
 Let $(V_g)_{g\in G}$ and $(W_g)_{g\in G}$ denote the $(C,F)$-actions of $G$ associated
 with $\Cal V$ and $\Cal W$ respectively.
 Then the following are satisfied:
 \roster
 \item"(P1)" 
 $(\psi_{\widetilde{\boldsymbol J}, \boldsymbol J} V_g\psi_{\widetilde{\boldsymbol J}, \boldsymbol J}^{-1})_{g\in G}=(W_g)_{g\in G}$,
 where
the $\psi_{\widetilde{\boldsymbol J}, \boldsymbol J}$ is the $(\widetilde{\boldsymbol J}, \boldsymbol J)$-chain equivalence. 
\item"(P2)" the $(C,F)$-action associated with $\Cal W_{\boldsymbol B}$ is $(\widehat\rho_{\boldsymbol B}W_g\widehat\rho_{\boldsymbol W}^{\,-1})_{g\in G}$, where
$\widehat\rho_{\boldsymbol B}$ is the $\boldsymbol B$-reduction mapping;

\item"(P3)" the $(C,F)$-action associated with $\Cal V_{\boldsymbol A}$ is $(\widehat\rho_{\boldsymbol A}V_g\widehat\rho_{\boldsymbol A}^{\,-1})_{g\in G}$, where
$\widehat\rho_{\boldsymbol A}$ is the $\boldsymbol A$-reduction mapping;
 \item"(P4)"
 the $(C,F)$-action associated with $\Cal T_{\boldsymbol A,\boldsymbol k}$ is
 $(\rho_{\boldsymbol A}\iota_{\boldsymbol k}T_g
\iota_{\boldsymbol k}^{-1}\rho_{\boldsymbol A}^{-1})_{g\in G}$,
where $\iota_{\boldsymbol k}$ is  the $\boldsymbol k$-telescoping
mapping and $\rho_{\boldsymbol A}$ is the $\boldsymbol A$-reduction mapping;
 \item"(P5)"
 the $(C,F)$-action associated with $\widetilde{\Cal T}_{{\boldsymbol B},\boldsymbol l}$ is
 $(\rho_{\boldsymbol B}\iota_{\boldsymbol l}\widetilde T_g
\iota_{\boldsymbol l}^{-1}\rho_{\boldsymbol B}^{-1})_{g\in G}$,
where $\iota_{\boldsymbol l}$ is  the $\boldsymbol l$-telescoping
mapping and $\rho_{\boldsymbol B}$ is the $\boldsymbol B$-reduction mapping.
\endroster
Since 
$\Cal V_{\boldsymbol A}=\Cal T_{\boldsymbol A,\boldsymbol k}$ and $\Cal W_{\boldsymbol B}=\widetilde{ \Cal T}_{\boldsymbol B,\boldsymbol l}$ it follows  from (P2)--(P5) that
$$
\widehat\rho_{\boldsymbol A} V_g\widehat\rho_{\boldsymbol A}^{-1}=\rho_{\boldsymbol A}\iota_{\boldsymbol k}T_g\iota_{\boldsymbol k}^{-1}\rho_{\boldsymbol A}^{-1}\quad\text{and}\quad
\widehat\rho_{\boldsymbol B} W_g\widehat\rho_{\boldsymbol B}^{-1}=\rho_{\boldsymbol B}\iota_{\boldsymbol l}\widetilde T_g\iota_{\boldsymbol l}^{-1}\rho_{\boldsymbol B}^{-1}.
\tag2-10
$$

 \proclaim{Theorem 2.2} Under the above notation, 
 $
 \phi=\iota_{\boldsymbol l}^{-1}\rho_{\boldsymbol B}^{-1}\widehat{ \rho}_{\boldsymbol B}\psi_{\widetilde{\boldsymbol J}, \boldsymbol J}
 \widehat{ \rho}_{\boldsymbol A}^{-1}\rho_{\boldsymbol A}\iota_{\boldsymbol k}.
 $
  \endproclaim

 \demo{Proof} Let 
 $
\theta:= \iota_{\boldsymbol l}^{-1}\rho_{\boldsymbol B}^{-1}\widehat{ \rho}_{\boldsymbol B}\psi_{\widetilde{\boldsymbol J}, \boldsymbol J}
 \widehat{ \rho}_{\boldsymbol A}^{-1}\rho_{\boldsymbol A}\iota_{\boldsymbol k}$.
 Then $\theta$ is a measurable isomorphism of $(X,\goth B,\mu)$ onto 
 $(\widetilde X,\widetilde{\goth B}, \widetilde\mu)$.
It follows from (P1) and \thetag{2-10} that $\theta T_g\theta^{-1}=\widetilde T_g$ for each $g\in G$.
 Take $n\in\Bbb N$.
 Then\footnote{We utilize \thetag{1-12} in the third equation.}
 $$
 \align
 \iota_{\boldsymbol k}([1_G]_{k_n})&=[1_G]_n,\\
 \widehat{ \rho}_{\boldsymbol A}^{\,-1}\rho_{\boldsymbol A}([1_G]_n)&=[1_G]_n,\\
\psi_{\widetilde{\boldsymbol J}, \boldsymbol J}([1_G]_n)&=[\widetilde J_{n}]_n,\\
 \widehat{ \rho}_{\boldsymbol B}\rho_{\boldsymbol B}^{-1}([\widetilde J_n]_n)&=[\widetilde J_n]_n\quad\text{and}\\
 \iota_{\boldsymbol l}^{-1}([\widetilde J_n]_n)&=[\widetilde J_{n}]_{l_n}.
 \endalign
 $$
Thus,  $\theta([1_G]_{k_n})=[\widetilde J_{n}]_{l_n}$.
On the other hand, \thetag{2-1} yields that
$$
\widetilde\mu(\phi([1_G]_{k_n})\triangle [\widetilde J_n]_{l_n})<\epsilon_n\mu([1_G]_{k_n}).
$$
Hence,
$$
\lim_{n\to\infty}\frac{\widetilde\mu(\theta([1_G]_{k_n})\triangle \phi([1_G]_{k_n}))}{\mu([1_G]_{k_n})} =0.
$$
 We now deduce from Proposition~1.9 (passing first to the $\boldsymbol k$-telescoping in $\Cal T$ and $\widetilde{\Cal T}$) that
 $\theta=\phi$ almost everywhere.
 \qed
 \enddemo

 \comment

 Let $\tau$ be a measure preserving mapping from a $(C,F)$-space $(X,\mu)$
 to another $(C,F)$-space $(\widetilde X,\widetilde\mu)$.
 Given $n,m>0$, we let $Y_{n,m}:=X_n\cap\phi^{-1}(\widetilde X_m)$.
 Then each $x$ belonging to $Y_{n,m}$ has
 infinitely many coordinates $(x_n,x_{n+1}, \dots)\in X_n$ and
 $\tau(x)$ has infinitely many coordinates $(\tau_m(x),\tau_{m+1}(x),\dots)\in \widetilde X_m$.

 \definition{Definition 2.3} 
  If for  each $n>0$, $m>0$ and $l\ge m$, the mapping
 $Y_{n,m}\ni x\mapsto\tau_l(x)$ depends, in fact, on finitely  many coordinates of $x$ almost everywhere then we say that $\tau$ is {\it finitary}.
  \enddefinition
 
 If $\tau$ is continuous then $\tau$ is finitary.
 Hence, all calibrating, telescoping and congruency mappings (and their inverses) are finitary.
 The reduction mapings are also  finitary.
If $\tau$ is a finitary measure preserving isomorphism from $X$ to $\widetilde X$
and $\tau'$ is a finitary
 measure preserving isomorphism from $\widetilde X$ to $ X'$  then
$\tau'\circ\tau$ 
is a finitary
 measure preserving isomorphism from $ X$ to $ X'$.
 It follows from Theorem~2.1 and~2.2 that if two $(C,F)$-actions are
 isomorphic then the corresponding isomorphism is finitary.
 
 \endcomment

 We illustrate Theorems~2.1 and 2.2 with the following example.
 
 \example {Example 2.3}
 Let $G=\Bbb Z$ and let
$\Cal T=(C_n,F_{n-1})_{n=1}^\infty$ be a $(C,F)$-sequence satisfying \thetag{1-1} and  \thetag{1-3}.
We assume that $F_n=\{0,1,\dots, h_n-1\}$, $n\in\Bbb N$, for an increasing sequence $(h_n)_{n=1}^\infty$ 
of positive integers.
Since $(F_n)_{n=1}^\infty$ is a F{\o}lner sequence in $\Bbb Z$, it follows that \thetag{1-5} holds.
We also assume that there is a sequence $(\epsilon_n)_{n=1}^\infty$ of positive reals 
and a sequence $(\beta_n)_{n=1}^\infty$ of positive integers such that
$\sum_{n=1}^\infty\epsilon_n<\infty$ and 
$$
\#((C_{n} +\beta_n)\cap C_n)>(1-\epsilon_n)\# C_n\quad\text{for each $n$.}\tag2-11
$$
 Then, of course, $\# C_n\to\infty$.
Denote by $T=(T_n)_{n
\in \Bbb Z}$ the 
 $(C,F)$-action associated with $\Cal T$.
 Let $(X,\mu)$ be the space of $T$.
  We now let 
  $$
  \alpha_n:=\beta_1+\dots+\beta_n.
  $$
 It is straightforward to verify that  
 $$
 \#((F_{n}+\alpha_n)\cap F_n)\ge (1-2\epsilon_n)\# F_{n}\quad\text{for each $n$.}\tag2-12
 $$
For $n>0$, take $x=(f_n,c_{n+1},\dots)\in X_n$  such that  $f_n\in (F_{n}-\alpha_n)\cap F_n$ and 
$c_{j}\in (C_{j} -\beta_{j})\cap C_{j}$, \dots
for all $j>n$.
We now set 
$$
\theta x:= (\alpha_n+f_n, c_{n+1}+\beta_{n+1}, c_{n+2}+\beta_{n+2}, \dots )\in X_n.
$$  
It follows from \thetag{2-11}, \thetag{2-12} and the Borel-Cantelli lemma
 that $\theta$ is a well-defined (mod 0)  measure  preserving invertible transformation 
 of $(X,\mu)$.
 Of course, $\theta\in C(T)$, i.e. $\theta$ is an isomorphism of $T$ with $T$.
 It is straightforward to verify that $T_{\alpha_n}\to\theta$ weakly.
 The latter means that $\lim_{n\to\infty}\mu(T_{\alpha_n} F\cap E)=\mu(\theta F\cap E)$
 for all subsets $E,F\subset X$.
 
Our purpose is to decompose $\theta$ into a product  of seven  ``elementary'' mappings as in Theorem~2.2.
Let
$$
\gather
k_n:=2n, \ l_n:=2n+1,\\
J_n:=(C_{2n}+\alpha_{2n})\cap F_{2n}\quad\text{and}\quad \widetilde J_n:=(C_{2n+1}-\alpha_{2n})\cap F_{2n+1}.
\endgather
$$
Then $J_n\subset F_{k_n}$ and  $\widetilde J_n\subset F_{l_n}$.
It is a routine to verify that  (i)--(iii) and (i)$'$--(iii)$'$ from the statement of Theorem~2.1
 hold for the sequence $(k_n,l_n, J_n,\widetilde J_n)_{n}$.
 \comment

We see that $C_{l_{n-1}+1,l_n}=C_{2n}+C_{2n+1}$.
Therefore, 
$$
\#(C_{l_{n-1}+1,l_n}\triangle (J_n+\widetilde J_n))\le 2\#((C_{2n}+\alpha_{2n})\setminus F_{2n})
+2\#((C_{2n+1}-\alpha_{2n})\setminus F_{2n+1}).
$$
This inequality and \thetag{2-13} imply that (iv)$'$ from the statement of Theorem~2.1
 holds for $(k_n,l_n, J_n,\widetilde J_n)_{n}$.
 On the other hand, $C_{k_n+1,k_{n+1}}=C_{2n+1}+C_{2n+2}$.
Since
 $$
\align
  \Bigg|\frac{\#((\widetilde J_n+J_{n+1})\triangle C_{k_n+1,k_{n+1}})}
 {\#C_{k_n+1,k_{n+1}}}&-
 \frac{\#((C_{2n+1}-\alpha_{2n}+C_{2n+2}+\alpha_{2n+2})\triangle C_{k_n+1,k_{n+1}})}
  {\#C_{k_n+1,k_{n+1}}}\Bigg|\\
  &\le 4\epsilon_n
  \endalign
 $$
and 
$$
C_{2n+1}-\alpha_{2n}+C_{2n+2}+\alpha_{2n+2}=C_{2n+1}+\beta_{2n+1}+C_{2n+2}+\beta_{2n+2},
$$  it follows from \thetag{2-12} that (iv) from the statement of  Theorem~2.1 holds.
Thus, all the conditions of Theorem~2.1 are satisfied for the sequence
$(k_n,l_n, J_n,\widetilde J_n)_{n}$.

\endcomment
We leave this verification to the reader.
Then, by Theorem~2.1, an isomorphism $\phi\in C(T)$ is well defined by
the sequence $(k_n,l_n, J_n,\widetilde J_n)_{n}$.
According to Theorem~2.2, $\phi$ is a composition of 7 elementary mappings:
$$
\phi=\iota_{\boldsymbol l}^{-1}\rho_{\boldsymbol B}^{-1}\widehat{ \rho}_{\boldsymbol B}\psi_{\widetilde{\boldsymbol J}, \boldsymbol J}
 \widehat{ \rho}_{\boldsymbol A}^{-1}\rho_{\boldsymbol A}\iota_{\boldsymbol k}
$$
that were introduced above the statement of Theorem~2.2.
Given $x\in X$, we now compute $\phi(x)$ coordinatewise.
Since the reduction mappings $\rho_{\boldsymbol B}$, $\widehat{ \rho}_{\boldsymbol B}$,
 $\widehat{ \rho}_{\boldsymbol A}$ and $\rho_{\boldsymbol A}$
 do not change coordinates of points from their domains,
we have to compute indeed only ``the actions'' of $\iota_{\boldsymbol k}$, 
$\psi_{\widetilde{\boldsymbol J}, \boldsymbol J}$ and $\iota_{\boldsymbol l}^{-1}$.
In view of  \thetag{2-11} and \thetag{2-12}, we can assume without loss of generality (i.e. dropping to a $\mu$-conull subset) that there is  $n=n(x)>0$ such that 
$$
x=(f_{2n-1},c_{2n}, c_{2n+1},\dots)\in X_{2n-1}\quad\text{ and $c_j+\beta_j, c_j-\beta_j\in C_j$ for each $j\ge 2n$.}
$$
This implies that $c_{2j}+\alpha_{2j}\in J_j$ and $c_{2j+1}-\alpha_{2j+1}\in \widetilde J_{j+1}$ for each $j\ge n$. 
Hence,
for each $j\ge n$,
$$
c_{2j}+c_{2j+1}=(c_{2j}+\alpha_{2j})+(c_{2j+1}-\alpha_{2j})\in J_j+\widetilde J_{j+1}.
$$
Since  
$
\iota_{\boldsymbol k}(x)=(f_{2n-1}, c_{2n}+c_{2n+1},
c_{2n+2}+c_{2n+3},\dots),
$
we obtain that
$$
\multline
\psi_{\widetilde{\boldsymbol J}, \boldsymbol J}
 \widehat{ \rho}_{\boldsymbol A}^{\,-1}\rho_{\boldsymbol A}\iota_{\boldsymbol k}(x)
\\
=\psi_{\widetilde{\boldsymbol J}, \boldsymbol J}
(f_{2n-1}, (c_{2n}+\alpha_{2n})+(c_{2n+1}-\alpha_{2n}),
(c_{2n+2}+\alpha_{2n+2})+(c_{2n+3}-\alpha_{2n+2}),\dots)\\
=
(f_{2n-1}+(c_{2n}+\alpha_{2n}),   (c_{2n+1}-\alpha_{2n})+
(c_{2n+2}+\alpha_{2n+2}), (c_{2n+3}-\alpha_{2n+2})+
(c_{2n+4}+\alpha_{2n+4}),\dots)\\
=
(f_{2n}+\alpha_{2n},   c_{2n+1}+\beta_{2n+1}+
c_{2n+2}+\beta_{2n+2},
 c_{2n+3}+\beta_{2n+3}
 +c_{2n+4}+\beta_{2n+4},\dots)
 \endmultline
 $$
 and 
 $$
 \align
 & 
 \iota_{\boldsymbol l}^{-1}\rho_{\boldsymbol B}^{-1}\widehat{ \rho}_{\boldsymbol B}(\psi_{\widetilde{\boldsymbol J}, \boldsymbol J}
 \widehat{ \rho}_{\boldsymbol A}^{\,-1}\rho_{\boldsymbol A}\iota_{\boldsymbol k}(x))
\\
 &=
(f_{2n}+\alpha_{2n},   c_{2n+1}+\beta_{2n+1},
c_{2n+2}+\beta_{2n+2},
 c_{2n+3}+\beta_{2n+3},
 c_{2n+4}+\beta_{2n+4},\dots).
\endalign
$$
It follows that $\phi=\theta$ almost everywhere.
Thus, $\theta=\iota_{\boldsymbol l}^{-1}\rho_{\boldsymbol B}^{-1}\widehat{ \rho}_{\boldsymbol B}\psi_{\widetilde{\boldsymbol J}, \boldsymbol J}
 \widehat{ \rho}_{\boldsymbol A}^{-1}\rho_{\boldsymbol A}\iota_{\boldsymbol k}$, as desired.

 \endexample

 \head 3. Factors of rank-one actions
 \endhead
 
 \subhead 3.1. Continuous proper factors of topological $(C,F)$-actions
 \endsubhead
 Let  $\Cal T=(C_n, F_{n-1})_{n>0}$ and $\widetilde{\Cal  T}=(\widetilde C_n, \widetilde F_{n-1})_{n>0}$
be two $(C,F)$-sequences satisfying 
\thetag{1-1} %\thetag{1-3} 
and \thetag{1-4}.
Denote by $T=(T_g)_{g\in G}$ and $\widetilde T=(\widetilde T_g)_{g\in G}$ the topological $(C,F)$-actions associated with $\Cal T$ and $\widetilde{\Cal T}$ respectively.
Let $X$ and $\widetilde X$ be the locally compact Cantor $(C,F)$-spaces  on which $T$ and $\widetilde T$ are determined respectively.

\definition{Definition 3.1} We say that $\widetilde {\Cal T}$ is a {\it quotient} of ${\Cal T}$ if
there is a sequence $\boldsymbol A:=(A_n)_{n=1}^\infty$ of finite subsets $A_n$ in $G$
such that  the following holds for each $n\ge 1$:
$$
\gather
F_{n-1}C_n\subset\widetilde F_nA_n\subset F_n,\tag3-1\\
\widetilde F_n^{-1}\widetilde F_n\cap A_n A_n^{-1}=\{1_G\},\tag3-2\\
A_nC_{n+1}=\widetilde C_{n+1}A_{n+1}\tag3-3
\endgather
$$
\enddefinition
 
 We now define a mapping $q_{\boldsymbol A}:X\to \widetilde X$.
 Let $x\in X$.
 Then there is $n\ge 0$ such that
  $x=(f_n,c_{n+1},c_{n+2},\dots)\in X_n\subset X$.
  Our purpose is to define an element $q_{\boldsymbol A}(x)$.
  %:=(\widetilde f_{n+1},\widetilde c_{n+2},\widetilde c_{n+3},\dots)\in\widetilde X_{n+1}$.
It follows from \thetag{3-1} that $f_nc_{n+1}\in \widetilde F_{n+1}A_{n+1}$.
In view of \thetag{3-2}, there exist a unique $\widetilde f_{n+1}\in\widetilde F_{n+1}$
 and a unique $a_{n+1}\in A_{n+1}$ such that 
 $$
 f_nc_{n+1}=\widetilde f_{n+1}a_{n+1}.
 $$
It follows from this and  \thetag{3-3} that  
$$
 f_nc_{n+1}c_{n+2}=\widetilde f_{n+1}a_{n+1}c_{n+2}=
 \widetilde f_{n+1}\widetilde c_{n+2}a_{n+2}
 $$
 for some  $\widetilde c_{n+2}\in\widetilde C_{n+2}$ and $a_{n+2}\in A_{n+2}$.
According to \thetag{3-2}, the elements $\widetilde c_{n+2}$ and 
$a_{n+2}$ are defined uniquely.
Continuing this procedure infinitely many times, we construct a sequence
$(\widetilde c_m)_{m>n+1}$ with ${\widetilde c}_m\in \widetilde C_m$ for each $m>n+1$.
We now set 
$$
q_{\boldsymbol A}(x):=(\widetilde f_{n+1},\widetilde c_{n+2},\widetilde c_{n+3},\dots)\in\widetilde X_{n+1}\subset\widetilde X
$$
It is a routine to verify that $q_{\boldsymbol A}$ is well defined as a mapping of $X$ to $\widetilde X$.
Of course, $q_{\boldsymbol A}$ is continuous
and  $q_{\boldsymbol A}(X_n)\subset \widetilde X_{n+1}$ for each  $n$.

We now show that $q_{\boldsymbol A}$ is onto. 
For that, it is sufficient to prove that $q_{\boldsymbol A}(X_{n+1})=\widetilde X_{n+1}$ for each $n>0$. 
Take a point $(\widetilde f_{n+1},\widetilde c_{n+2},\widetilde c_{n+3},\dots)\in\widetilde X_{n+1}$.
For each $m>n+1$ and an element $a_m\in A_m$, we apply \thetag{3-3} repeatedly and 
then \thetag{3-1} to determine uniquely  the following elements:   $a_{m-1}\in A_{m-1}$, \dots, $a_{n+1}\in A_{n+1}$,  $c_m\in C_m$, \dots, $c_{n+1}\in C_{n+1}$ and $f_n\in F_n$ such that
$$
\align
\widetilde f_{n+1}\widetilde c_{n+2}\cdots\widetilde c_{m} a_m &=
\widetilde f_{n+1}\widetilde c_{n+2}\cdots \widetilde c_{m-1}a_{m-1} c_{m} \\
&=
\widetilde f_{n+1}\widetilde c_{n+2}\cdots  \widetilde c_{m-2}a_{m-2} c_{m-1} c_{m} \\
&\cdots \\ &=\widetilde f_{n+1}a_{n+1}c_{n+2}\cdots c_m\\
&=f_{n+1}c_{n+2}\cdots c_m.
\endalign
$$
 Thus, to each $a_m\in A_m$ we put in correspondence a finite sequence $(a_j)_{j=n+1}^m$ such that
 $a_j\in A_j$, $a_{j-1}^{-1}\widetilde c_j a_j\in C_j$ for $j=n+2,\dots,m$ and $\widetilde f_{n+1}a_{n+1}\in F_{n+1}$. 
 Since $F_{n+1}$ and $A_j$ is finite for each $j$ and $\bigcup_{m>n+1}A_m$ is infinite, it follows that there exists an infinite sequence $(a_j)_{j=n+1}^\infty$
 such that $f_{n+1}:=\widetilde f_{n+1}a_{n+1}\in F_{n+1}$, 
 $a_j\in A_j$ and $c_j:=a_{j-1}^{-1}\widetilde c_j a_j\in C_j$ for each $j>n+1$. 
 Then $x:=(f_{n+1},c_{n+2}, c_{n+3},\dots)$ belongs to $X_{n+1}$ and $q_{\boldsymbol A}(x)=\widetilde x$, as desired.
 
 We now prove that $q_{\boldsymbol A}T_g=\widetilde T_gq_{\boldsymbol A}$ for each $g\in G$.
 Let $x\in X$ and $g\in G$.
 Since~\thetag{1-4} holds for $\Cal T$ and $\widetilde{ \Cal T}$, there is $n>0$ such that 
 $$
 \align
 x&=(f_n,c_{n+1},\dots)\in X_n,\quad
  f_n,gf_n\in F_n,\\
 q_{\boldsymbol A} x&=(\widetilde f_{n+1},\widetilde c_{n+2},\dots)\in\widetilde X_{n+1}\quad \text{and}\quad
 \widetilde f_{n+1},g\widetilde f_{n+1}\in \widetilde F_{n+1}.
 \endalign
 $$
 It follows from the definition of $q_{\boldsymbol A}$ that
 $
 f_nc_{n+1}=\widetilde f_{n+1}a_{n+1}  $
 for some $a_{n+1}\in A_{n+1}$.
 Hence
 $
 gf_nc_{n+1}=g\widetilde f_{n+1}a_{n+1}  $.
 This yields that 
 $$
 q_{\boldsymbol A}(T_gx)=q_{\boldsymbol A}(gf_n,c_{n+1},\dots) =
 (g\widetilde f_{n+1},\widetilde c_{n+2},\dots)=\widetilde T_gq_{\boldsymbol A}(x),
 $$
as desired.

\definition{Definition 3.2} We call $q_{\boldsymbol A}$ the {\it $\boldsymbol A$-quotient mapping.}
\enddefinition

 It is straightforward to verify that 
 $$
 q_{\boldsymbol A}^{-1}([\widetilde f]_n)=[\widetilde f A_n]_n\quad\text{for each $\widetilde f\in\widetilde F_{n}$ and $n>0$.}\tag3-4
 $$
 Hence, the $q_{\boldsymbol A}$-inverse image of each compact open  subset in $\widetilde X$ is compact.
  Since the compact open subsets are a base of the topology in $\widetilde X$, if follows that the $q_{\boldsymbol A}$-inverse image of each compact  subset in $\widetilde X$ is compact in $X$.
  Hence,   $q_{\boldsymbol A}$ is proper.
  
  Let $\mu$ and $\widetilde\mu$ be the  Haar measures on $X$ and $\widetilde X$ respectively.
   %(i.e. we assume that $\mu(X)=\widetilde \mu(\widetilde X)=1$).
   Since $q_{\boldsymbol A}$ is proper, the measure $\mu\circ q_{\boldsymbol A}^{-1}$ is Radon.
Since $\mu$ is invariant under $T$ and $q_{\boldsymbol A}$ is equivariant, it follows that $\mu\circ q_{\boldsymbol A}^{-1}$ is invariant under $\widetilde T$.
Since $\widetilde T$ is Radon uniquely ergodic, we obtain that 
$\mu\circ q_{\boldsymbol A}^{-1}= d\widetilde\mu$ for some constant $d>0$.

 In the following theorem we find necessary and sufficient conditions under which  a
continuous $(C,F)$-action on a locally compact Cantor space is a proper continuous factor of another
 continuous $(C,F)$-action on a locally compact Cantor space.
 These conditions are given  in terms of the underlying $(C,F)$-parameters.
 Moreover, an explicit formula for the factor mappings is obtained.
 
 \proclaim{Theorem 3.3}
 Let  $\Cal T=(C_n, F_{n-1})_{n>0}$ and $\widetilde{\Cal  T}=(\widetilde C_n, \widetilde F_{n-1})_{n>0}$
be two $(C,F)$-sequences satisfying 
\thetag{1-1},  \thetag{1-4} and \thetag{1-7}.
Denote by $T=(T_g)_{g\in G}$ and $\widetilde T=(\widetilde T_g)_{g\in G}$ the topological $(C,F)$-actions associated with $\Cal T$ and $\widetilde{\Cal T}$ respectively.
Let $X$ and $\widetilde X$ be the locally compact Cantor $(C,F)$-spaces  on which $T$ and $\widetilde T$ are determined respectively.
A proper continuous onto mapping $\theta:X\to \widetilde X$ that intertwines $T$ with $\widetilde T$ exists
if and only if  there are
 an increasing sequence of integers $\boldsymbol k=(k_n)_{n=0}^\infty$  with $k_0=0$ and
 a sequence $\boldsymbol A=(A_n)_{n=1}^\infty$ of finite subsets in $G$
such that  \thetag{3-1}--\thetag{3-3}  are satisfied with
the $\boldsymbol k$-telescoping of $\Cal T$ in place of $\Cal T$.
Moreover,
$\theta=q_{\boldsymbol A}\iota_{\boldsymbol k}$, where $\iota_{\boldsymbol k}$ is the $\boldsymbol k$-telescoping mapping and $q_{\boldsymbol A}$ is the $\boldsymbol A$-quotient mapping.
 \endproclaim

 \demo{Proof}
 Since the ``if'' part of the statement of the theorem has been proved above (at the beginning of \S3.1), it remains  to prove the  ``only if'' part.
  Thus, $\theta$ is given.
  Our goal is 
to construct $(A_n)_{n=1}^\infty$ and $(k_n)_{n=1}^\infty$
satisfying the required conditions. 
We will do this
 inductively.

 Since $\theta$ is continuous and proper, for each $n\ge 0$, the subset $\theta^{-1}[1_G]_n$ is compact and open.
 Hence, on the $n$-th step,  we can choose $k_n>k_{n-1}$ and  a subset $A_n\subset F_{k_n}$ such that 
 $\theta^{-1}[1_G]_n=[A_n]_{k_n}$.
 Of course, 
 $$
 \theta^{-1}\widetilde X_n=\bigsqcup_{g\in \widetilde F_n}\theta^{-1}(\widetilde T_g[1_G]_n)=\bigsqcup_{g\in \widetilde F_n}T_g[A_n]_{k_n}.
 $$
 Since \thetag{1-4} holds, we can assume without loss of generality (increasing $k_n$ if necessary)
  that $\widetilde F_nA_n\subset F_{k_n}$.
 Moreover,
 $$
 gA_n\cap g'A_n=\emptyset\quad\text{for all $g,g'\in\widetilde F_n$, $g\ne g'$.}\tag3-5
 $$
 Since $X=\bigcup_{n>0}\theta^{-1}(\widetilde X_n)$, we can assume additionally that
 $\theta^{-1}(\widetilde X_n)\supset X_{k_{n-1}}$.
 Thus, 
 $$
 F_{k_{n-1}}C_{k_{n-1}+1,k_n}\subset
 \widetilde F_nA_n\subset F_{k_n}.
 \tag3-6
 $$
 Also,  $\theta^{-1}[1_G]_n=\theta^{-1}[\widetilde C_{n+1}]_{n+1}$.
 Hence 
 $$
 A_{n}C_{k_n+1,k_{n+1}}=\widetilde C_{n+1}A_{n+1}.
 \tag3-7
 $$
 We now let $\boldsymbol k:=(k_n)_{n\ge 0}$ and $\boldsymbol A:=(A_n)_{n=1}^\infty$.
Denote by $\Cal T'$ the $\boldsymbol k$-telescoping of $\Cal T$.
Then $\widetilde{ \Cal T}$ is a quotient of $\Cal T'$
and \thetag{3-5}, \thetag{3-6} and \thetag{3-7} are analogues of \thetag{3-2}, \thetag{3-1} and
\thetag{3-3} respectively.
Therefore,  the mapping $q_{\boldsymbol A}\iota_{\boldsymbol k}$ is equivariant, i.e. $q_{\boldsymbol A}\iota_{\boldsymbol k}T_g=\widetilde T_gq_{\boldsymbol A}\iota_{\boldsymbol k}$
for each $g\in G$.
In view of \thetag{3-4},
$$
(q_{\boldsymbol A}\iota_{\boldsymbol k})^{-1}[1_G]_n=[A_n]_{k_n}=\theta^{-1}[1_G]_n\quad\text{for each $n>0$}.
$$
Since $q_{\boldsymbol A}\iota_{\boldsymbol k}$ and $\theta$ are both equivariant, it follows that 
$$
(q_{\boldsymbol A}\iota_{\boldsymbol k})^{-1}O=\theta^{-1}O
\quad\text{for each cylinder $O\subset \widetilde X$.}
$$
It follows that  $q_{\boldsymbol A}\iota_{\boldsymbol k}=\theta$.
 \qed
 \enddemo

 \subhead 3.2. Measurable  factors of measure theoretical $(C,F)$-actions
 \endsubhead
 Let  $\Cal T=(C_n, F_{n-1})_{n>0}$ and $\widetilde{\Cal  T}=(\widetilde C_n, \widetilde F_{n-1})_{n>0}$
be two $(C,F)$-sequences satisfying 
\thetag{1-1}, \thetag{1-3}, \thetag{1-5} and \thetag{1-7}.
Denote by $T=(T_g)_{g\in G}$ and $\widetilde T=(\widetilde T_g)_{g\in G}$ the $(C,F)$-actions associated with $\Cal T$ and $\widetilde{\Cal T}$ respectively.
%Then $T$ and $\widetilde T$ are finite measure preserving.
Let $(X,\goth B,\mu)$ and $(\widetilde X,\widetilde{\goth B},\widetilde\mu)$ be the corresponding standard probability $(C,F)$-spaces  on which $T$ and $\widetilde T$ are determined.
Thus, $\mu$ and $\widetilde\mu$ are the normalized Haar measures for the tail equivalence relations on $X$ and $\widetilde X$ respectively.
Fix a decreasing sequence $(\epsilon_n)_{n=1}^\infty$ of positive reals such that  $\epsilon_1>1$ and $\sum_{n=2}^\infty\epsilon_n<0.2$.
Without loss of generality (passing to a telescoping $\widetilde T$, if necessary) we may assume 
that
$$
\widetilde\mu(\widetilde X_{n})>1-\frac{\epsilon_n}2.\tag3-8
$$
The following theorem provides necessary and sufficient conditions under which $\widetilde T$
is a measure theoretical factor of $T$.
The conditions are given in terms of the underlying $(C,F)$-parameters.

%Denote by $T=(T_g)_{g\in G}$ and $\widetilde T=(\widetilde T_g)_{g\in G}$ the  $(C,F)$-actions associated with $\Cal T$ and $\widetilde{\Cal T}$ respectively.
%Let $(X,\mu)$ and $(\widetilde X,\widetilde\mu)$ be the probability $(C,F)$-spaces  on which $T$ and $\widetilde T$ are defined.

 \proclaim{Theorem 3.4}  $\widetilde T$ is isomorphic to a  (measure theoretical) factor of $T$ if and only if there exist an increasing sequence
$
0=k_0<k_1<k_2<\cdots
$
 of non-negative integers and subsets $J_n\subset F_{k_{n}}$ such that
\roster
\item"(i)"
$\widetilde F_{n}  J_n \subset  F_{k_{n}}$,
\item"(ii)" 
the mapping $\widetilde F_{n}\times J_n \ni (\widetilde f,f) \mapsto \widetilde ff \in F_{k_{n}}$ is one-to-one, 
%(i.e. $(F_{k_{n+1}} - F_{k_{n+1}}) \cap (\widetilde J_n - \widetilde J_n) = \{0\}$),
\item"(iii)" 
$\dfrac{\# F_{k_n} - \# \widetilde F_{{n}} \#  J_n}{\#  F_{k_n}} < \epsilon_n$ and
\item"(iv)" 
$\dfrac{\#((J_{n-1}  C_{k_{n-1}+1,k_{n}}) \Delta \widetilde C_{n }J_{n})}{\# \widetilde C_{n }\# J_{n}} < 2\epsilon_{n-1}$
\endroster
for each $n\ge 1$.
\endproclaim

\demo{Proof} We first prove the ``only if" claim.
 Let $\phi:X\to \widetilde X$ be a measure preserving isomorphism that intertwines
$T$ with $\widetilde T$.
We will construct the desired objects inductively.
On the first step we  let $ J_0:=\{1_G\}$.
Suppose that for some $n\in\Bbb N$, we have already constructed  integers $(k_j)_{j=0}^{n-1}$  and 
subsets $(J_m)_{m=0}^{n-1}$   that satisfy (i)--(iv).
%We will assume, in addition, that the following inequality holds:
%\roster
%\item"(v)"
% $\widetilde\mu(\widetilde X\setminus \widetilde X_{n-1}) < \epsilon_{n-1}$.
% \endroster
Our purpose is to find an integer  $k_{n}$  such that $k_{n}>k_{n-1}$ and  a subset $J_{n}\subset F_{k_n}$
  for which (i)--(iv)  are satisfied.
Given  $l>k_{n-1}$,
we let 
$$
 F_{l}^\circ:=\{f\in  F_{l}\mid \widetilde F_{n}f\subset  F_{l}\}.
 $$
 Since $( F_l)_{l=1}^\infty$ is a F{\o}lner sequence, $\#  F_{l}^\circ/\# F_{l}\to 1$ as $l\to\infty$.
The sequence of  rings (of cylinders) $\{[ A]_{l}\mid A\subset F_l^0\}$   approximates the entire $\sigma$-algebra ${\goth B}$ as $l\to\infty$.
Hence there is $k_{n}>k_{n-1}$ and a subset $J_{n}\subset F_{k_n}^\circ$ such that
$$
\align
\mu([J_n]_{k_n}\triangle \phi^{-1}([1_G]_n))&< \frac{\epsilon_n}2\mu(\phi^{-1}([1_G]_n)),\tag3-9 \\
\min_{f\in J_n}\mu([f]_{k_n}\cap \phi^{-1}([1_G]_n))&>0.5\mu(([f]_{k_n})\quad\text{and}\tag3-10
\endalign
$$
The inclusion $J_{n}\subset F_{k_n}^\circ$ implies  (i). 
It is a routine to show that \thetag{3-10} implies (ii): a cylinder $[\widetilde f f]_{k_n}$ is mostly filled
with $\phi^{-1}([\widetilde f]_n)$ and $\phi^{-1}([\widetilde f]_n)\cap \phi^{-1}([\widetilde f']_n)=\emptyset$ whenever $\widetilde f\ne\widetilde f'$.
We apply \thetag{3-9} to obtain the following:
$$
\align
\widetilde\mu(\widetilde X_{n})&=\mu(\phi^{-1}([\widetilde F_n]_n))\\
&=
\mu\Bigg(\bigsqcup_{\widetilde f\in\widetilde F_{n}}T_{\widetilde f}\phi^{-1}([1_G]_n)\Bigg)\\
&=
\# \widetilde F_n\mu(\phi^{-1}([1_G]_n)
)\\
&\le 
(1-0.5\epsilon_n)^{-1}
\# \widetilde F_n\mu([J_n]_{k_n})\\
&=(1-0.5\epsilon_n)^{-1}\# \widetilde F_n\# J_n\mu([1_G]_{k_n})\\
&<\frac{\# \widetilde F_n\# J_n}{(1-0.5\epsilon_n)\#F_{k_n}}.
\endalign
$$
Therefore, \thetag{3-8} entails that 
$$
1-\epsilon_n<(1-0.5\epsilon_n)^2<\frac{\# \widetilde F_n\# J_n}{\# F_{k_n}}.
$$
This yields (iii).

Since $[1_G]_{n-1}=[\widetilde C_{n}]_{n}$, it follows that 
$$
\phi^{-1}([1_G]_{n-1})=\bigsqcup_{c\in \widetilde C_{n}}T_c\phi^{-1}([1_G]_{n}).
$$
Therefore, applying \thetag{3-9} we obtain that
$$
\align
0&=\mu(\phi^{-1}([1_G]_{n-1})\triangle \bigsqcup_{c\in \widetilde C_{n}}T_c\phi^{-1}([1_G]_{n}))\\
&\ge
\mu\bigg([J_{n-1}]_{k_{n-1}}\triangle \bigsqcup_{c\in \widetilde C_{n}}T_c[J_{n}]_{k_{n}}\bigg)
- \frac{\epsilon_{n-1}}2\mu(\phi^{-1}([1_G]_{n-1}))-
 \frac{\epsilon_{n}\#\widetilde C_n}2\mu(\phi^{-1}([1_G]_{n}))\\
 &=
 \mu([J_{n-1}C_{k_{n-1}+1,k_n}]_{k_n} \triangle[\widetilde C_n J_n]_{k_n})-
 \frac{\epsilon_{n-1}+\epsilon_n}2\mu(\phi^{-1}([1_G]_{n-1}))\\
 &> \mu([(J_{n-1}C_{k_{n-1}+1,k_n}) \triangle (\widetilde C_n J_n)]_{k_n})-\epsilon_{n-1}(1+\epsilon_{n-1})\mu([J_{n-1}]_{k_{n-1}}).
\endalign
$$
Hence,  
$$
\frac{\#((J_{n-1}  C_{k_{n-1}+1,k_{n}}) \Delta \widetilde C_{n }J_{n})}{\#(J_{n-1}  C_{k_{n-1}+1,k_{n}})} < \epsilon_{n-1}(1+\epsilon_{n-1}).
$$
This inequality implies  (iv).
Thus, the ``only if'' claim is proved.

 We now prove the ``if'' claim.
 Thus, suppose that there exist two sequences $(k_n)_{n=0}^\infty$ and $(J_n)_{n=0}^\infty$
 such that (i)--(iv) are satisfied. 
 Let
 $$
 \multline
 Y_n:=\{(f, c)\in F_{k_n}\times C_{k_n+1,k_{n+1}}\mid f= \widetilde fj_n\text{ and } j_nc=\widetilde cj_{n+1}\\
 \text{ for some $\widetilde f\in \widetilde F_n,
  j_n\in J_n,\widetilde c\in\widetilde C_{n+1},j_{n+1}\in J_{n+1}$}\}
 \endmultline
 $$
 and let 
 $$
 Y_n^+:=\{(f_{k_n},c_{k_n+1},c_{k_n+2},\dots)\in X_{k_n}\mid (f_{k_n}, c_{k_n+1}\cdots c_{k_{n+1}})\in Y_n\}.
 $$
 It follows from (iii) and (iv) that
 $$
 \frac{\# Y_n}{\# F_{k_n}\# C_{k_n+1,k_{n+1}}}>1-3\epsilon_{n}.
 $$
Hence,
  $$
\sum_{n=1}^\infty3\epsilon_n>  \sum_{n=1}^\infty\bigg(1-\frac{\# Y_n}{\# F_{k_n}\# C_{k_n+1,k_{n+1}}}\bigg)
= \sum_{n=1}^\infty \frac{\mu(X_{k_n})-\mu(Y_n^+)}{\mu(X_{k_n})}.
  $$
This yields  that $\sum_{n=1}^\infty \mu(X_{k_n}\setminus Y_n^+)<\infty$.
 Hence,  it follows from the Borel-Cantelli that
 for a.e. $x\in X$, there exists $N>0$ such that $x=(f_{k_n},c_{k_n+1},c_{k_n+2},\dots)\in Y_n^+$ for each $n\ge N$.
This means that
 $$
 (f_{k_n}, c_{k_n+1}\cdots c_{k_{n+1}})\in Y_n,
 $$
i.e.  there are unique $\widetilde f_n\in\widetilde F_n$, $\widetilde c_{n+1}\in \widetilde C_{n+1}$
 and $j_{n+1}\in J_{n+1}$ such that
 $$
 f_{k_n} c_{k_n+1}\cdots c_{k_{n+1}}=\widetilde f_n\widetilde c_{n+1}j_{n+1}.\tag3-11
 $$
Since $x\in Y_{n+1}^+$, we also have that
 $$
  f_{k_n} c_{k_n+1}\cdots c_{k_{n+2}}=\widetilde f_{n+1}\widehat j_{n+1}c_{k_{n+1}+1}\cdots c_{k_{n+2}}
  \tag3-12
 $$
 for some $\widetilde f_{n+1}\in \widetilde F_{n+1}$ and $\widehat j_{n+1}\in J_{n+1}$.
 It follows from \thetag{3-11} and \thetag{3-12} that
 $$
 \widetilde f_n\widetilde c_{n+1}j_{n+1}=\widetilde f_{n+1}\widehat j_{n+1}.
 $$
Using (ii) we obtain that $\widetilde f_n\widetilde c_{n+1}=\widetilde f_{n+1}$.
 This equality for each $n\ge N$.
Therefore, 
 $$
 \theta x:=(\widetilde f_N,\widetilde c_{N+1}, \widetilde c_{N+2},\dots).
 $$
 is well defined as point of $\widetilde X$.
 Of course,
$\theta$ is a Borel mapping from (a conull subset of) $X$ to $\widetilde X$.
 It is straightforward to check that $\theta$ is equivariant: $\theta T_g=\widetilde T_g\theta$
 for each $g\in G$. 
 Hence the probability measure $\mu\circ\theta^{-1}$ on $\widetilde X$ and invariant under $\widetilde T$.
Since $\widetilde\mu$ is finite, $\widetilde T$ is uniquely ergodic.
Hence, $\mu\circ\theta^{-1}=\widetilde\mu$.
In particular, $\theta$ is onto (mod 0).
Thus, we have proved that $(\widetilde X,\widetilde \mu,\widetilde T)$ is a factor of $(X,\mu,T)$.
 \qed
\enddemo
 
 \subhead 3.3. Odometer factors
 \endsubhead
We will show that if  $G=\Bbb Z$  and $\widetilde T$ is an odometer,
then the statement  of Theorem~3.4 is equivalent to  the description of  odometer factors of rank-one maps 
from  \cite{Fo--We}  and \cite{DaVi}.

Let $(d_n)_{n=0}^\infty$ be a sequence of integers  such that $d_0=1$ and $d_n\ge 2$ if $n>0$.
We set for each $n\ge 0$,
$$
\widetilde C_{n+1}:=\{d_0\cdots d_{n}j\mid 0\le j< d_{n+1}\}
\quad\text{and}\quad 
\widetilde F_n:=\{0,1,\dots, d_0\cdots d_n-1\}.
$$
Then the sequence $(\widetilde C_{n+1},\widetilde F_n)_{n=0}^\infty$
satisfies  \thetag{1-1}, \thetag{1-3}, \thetag{1-5} and \thetag{1-7}.
Denote by $\widetilde T$ the corresponding $(C,F)$-action of $\Bbb Z$.
Of course, $\widetilde T$ is an odometer  and the discrete spectrum of
$\widetilde T$ is  
$\big\{e^{\frac{2\pi im}{d_1\cdots d_l}}\mid m\in\Bbb Z, l\in\Bbb N\big\}\subset\Bbb T$.
The following claim was first proved in \cite{Fo--We}  and then in \cite{DaVi} (in different, but equivalent terms).

\proclaim{Fact F}
Let $T$ be a  $(C,F)$-action of $\Bbb Z$ associated with a sequence 
$\Cal T=( C_{n+1}, F_n)_{n=0}^\infty$
satisfying~\thetag{1-1}, \thetag{1-3},
 \thetag{1-5} and \thetag{1-7}.
Then
$\widetilde T$ is a factor of $T$ if and only if there is an increasing
 sequence $0=k_0<k_1<k_2<\cdots$ of integers such that  
 %the  $(k_n)_{n=0}^\infty$-telescoping of $\Cal T$ satisfies  the condition
$$
\sum\limits_{n>0} \dfrac{\#\{c\in C_{k_n + 1, k_{n+1}} \mid c \neq 0 \mod{d_1d_2\ldots d_n}\}}{\# C_{k_n + 1, k_{n+1}}} < \infty.\tag3-13
$$
\endproclaim

We now show that Fact~F   is a corollary from   Theorem~3.4.
For that we will need some notation.
 Given two finite subsets $A,B\subset\Bbb Z$ and $\epsilon>0$, we write $A\approx_\epsilon B$ if
 $\#(A\triangle B)<\epsilon{\# B}$.
 It is straightforward to verify that 
 \roster
 \item"---"  if  $A\approx_\epsilon B$ then  if  $B\approx_{\epsilon/(1-\epsilon)} A$;
 \item"---" if  $A\approx_\epsilon B$ and $B\approx_\delta D$ then $A\approx_{\epsilon+\delta+\epsilon\delta} D$;
 \item"---"
 if $A\approx_\epsilon B$ and $C$ is another finite subset of $\Bbb Z$ such that $(A+c)\cap (A+c')=(B+c)\cap (B+c')=\emptyset$ whenever
 $c,c'\in C$ and $c\ne c'$ then $(A+C)\approx_\epsilon (B+C)$.
 \endroster
 We say that two subsets $A,B\subset\Bbb Z$ {\it do not overlap} if either $\max A<\min B$ or
 $\max B<\min A$.
 
 Since the ``if'' part of~Fact~F is straightforward (see \cite{DaVi} for details),
it remains to deduce  the ``only if" part of~Fact~F from  Theorem~3.4.
Thus, let $T$ the the $(C,F)$-action associated with $\Cal T$ and let
 let  $\widetilde T$ be a factor of $T$ as in Fact~F.
Fix a decreasing sequence $(\epsilon_n)_{n=1}^\infty$ of positive reals such that $1>\epsilon_n>6\sum_{j>n}\epsilon_j$ for each $n>0$.
Then, by Theorem~3.4,
there exist an increasing sequence
$
0=k_0<k_1<k_2<\cdots
$
 of non-negative integers and subsets $J_n\subset F_{k_{n}}$, $n\in\Bbb N$, such that (i)--(iv) of Theorem~3.4 hold.
 By Theorem~3.4(iv), for each $n>0$,
 $$
 J_{n-1} + C_{k_{n-1}+1,k_{n}}\approx_{2\epsilon_{n-1}} \widetilde C_{n }+J_{n}.
 $$
 Hence, 
 $$
 J_{n-1} + C_{k_{n-1}+1,k_{n}}+C_{k_{n}+1,k_{n+1}}\approx_{2\epsilon_{n-1}} \widetilde C_{n }+J_{n}+C_{k_{n}+1,k_{n+1}}.\tag3-14
 $$
 On the other hand, $J_{n}+C_{k_{n}+1,k_{n+1}}\approx_{2\epsilon_{n}} \widetilde C_{n +1}+J_{n+1}$
 and, hence,  
 $$
 \widetilde C_{n }+J_{n}+C_{k_{n}+1,k_{n+1}}\approx_{2\epsilon_{n}} \widetilde C_{n }+ \widetilde C_{n +1}+J_{n+1}.\tag3-15
 $$
 We deduce from \thetag{3-14} and \thetag{3-15} that
 $$
  J_{n-1} + C_{k_{n-1}+1,k_{n}}+C_{k_{n}+1,k_{n+1}}\approx_{2\epsilon_{n-1}+6\epsilon_n}
  \widetilde C_{n }+ \widetilde C_{n +1}+J_{n+1}.
 $$
Using this argument repeatedly, we obtain that for each $m>n$,
$$
  J_{n-1} + \sum_{j=n-1}^{m-1}C_{k_{j}+1,k_{j+1}}\approx_{2\epsilon_{n-1}
  +6\sum_{j=n}^{m-1}\epsilon_j}
 \Bigg(\sum_{j=n}^m \widetilde C_{j }\Bigg)+J_{m}.
 $$
Since $\epsilon_{n-1}>6\sum_{j=n}^{m-1}\epsilon_{j}$, it follows that
$$
  J_{n-1} + \sum_{j=n-1}^{m-1}C_{k_{j}+1,k_{j+1}}\approx_{3\epsilon_{n-1}}
 \Bigg(\sum_{j=n}^m \widetilde C_{j }\Bigg)+J_{m}.\tag3-16
 $$
Let 
$$
C:=J_{n-1}+\sum_{j=n}^{m-1}C_{k_{j}+1,k_{j+1}}\quad\text{and}\quad\widetilde C:=\sum_{j=n}^m \widetilde C_{j }.
$$
Then we can rewrite \thetag{3-16} as
$$
C_{k_{n-1}+1,k_n}+C\approx_{3\epsilon_{n-1}}\widetilde C+J_{m}.
$$
We note that 
$$
\widetilde C=\{jd_1\cdots d_n\mid j=0,1,\dots, (d_{n+1}\cdots d_m)-1\}
$$
is an arithmetic sequence with a common difference $d_1\cdots d_n$.
Let 
$$
L:=\max C_{k_{n-1}+1,k_n}
\quad\text{and}\quad
\widetilde C^\circ:=\{\widetilde c\in\widetilde C\mid L<\widetilde c< d_1\cdots d_m-L\}.
$$
We choose $m$ large so that $\widetilde C+J_{m}\approx_{\epsilon_{m}}\widetilde C^\circ+J_{m}$
and, hence,
$$
C_{k_{n-1}+1,k_n}+C\approx_{4\epsilon_{n-1}}\widetilde C^\circ+J_{m}.
$$
%We note that if $c+d=\widetilde c+j$ then $C_{k_{n-1}+1,k_n} +d\subset \widetilde C+j$.
Since the subsets  $(C_{k_{n-1}+1,k_n}+c)_{c\in C}$ are mutually disjoint, there exists  $c\in C$ such that 
$$
\#\big( (C_{k_{n-1}+1,k_n}+c)\cap (\widetilde C^\circ+J_{m})\big)>(1-4\epsilon_{n-1})\# C_{k_{n-1}+1,k_n}.\tag3-17
$$
Since
\roster
\item"---"
 the subsets $(\widetilde C+j)_{j\in J_{m}}$ do not pairwise overlap and 
 \item"---" the diameter of the set $C_{k_{n-1}+1,k_n}+c$ is $L$,
 \endroster
it follows that 
  there is  a unique $j\in J_{m}$ such that
$$
(C_{k_{n-1}+1,k_n}+c)\cap (\widetilde C^\circ+J_{m})\subset \widetilde C+j.\tag3-18
$$
From \thetag{3-17} and \thetag{3-18} we deduce that 
$$
\#\big( C_{k_{n-1}+1,k_n}\cap (\widetilde C+j-c)\big)>(1-4\epsilon_{n-1})\# C_{k_{n-1}+1,k_n}.
 $$
 Let $i_n:=j-c$.
An integer  $a$ belongs to $\widetilde C+j-c$
 if and only if $a-i_n$ is divisible by $d_1\cdots d_n$.
 Therefore,
 $$
\frac{
\# \{a\in C_{k_{n-1}+1,k_n}\mid a= i_n\mod d_1\cdots d_n\}
}
{\# C_{k_{n-1}+1,k_n}}
>1-4\epsilon_{n-1}.
 $$
 Thus, there is a sequence $(i_n)_{n=1}^\infty$ of integers such that
 $$
 \sum_{n=1}^
 \infty
 \frac{
\# \{a\in C_{k_{n-1}+1,k_n}\mid a\ne i_n\mod d_1\cdots d_n\}
}
{\# C_{k_{n-1}+1,k_n}}<\infty.
 $$
 Passing to a further telescoping of the $(k_n)_{n=0}^\infty$-telescoping of $\Cal T$, we can achieve that
 $i_1=i_2=\cdots=0$ (see \cite{DaVi} for details).
 Thus, \thetag{3-13} holds.
The ``only if" part of Fact~F
 is proved.

 \head 4. On classification of rank-one actions 
 \endhead

\subhead 4.1. Finite measure preserving rank-one actions
\endsubhead
Fix a countable discrete amenable group $G$ and a  standard probability space $(X,\mu):=([0,1], \text{Leb})$.
Denote by Aut$(X,\mu)$ the group of $\mu$-preserving transformations of $X$.
Endow Aut$(X,\mu)$ with the weak topology.
Then Aut$(X,\mu)$ is a Polish group.
Endow  the infinite product space $\text{Aut}(X,\mu)^G$  with the infinite product of the weak topologies 
 on $\text{Aut}(X,\mu)$.
 Denote by $\Cal A_G$ the set of  measure preserving $G$-actions on $(X,\mu)$. 
 Each element of $\Cal A_G$ is a homomorphism from $G$ to 
 $\text{Aut}(X,\mu)$.
Hence, $\Cal A_G$ is a subset of  $\text{Aut}(X,\mu)^G$.
It is straightforward to verify that this subset is closed.
 Hence, $\Cal A_G$ is a Polish space in the induced topology.
 The group  $\text{Aut}(X,\mu)$ acts on $\Cal A_G$ by conjugation.
This action is continuous.
Two $G$-actions from $\Cal A_G$ are isomorphic if and only if they
belong to the same  $\text{Aut}(X,\mu)$-orbit.
 
 Let $\goth F_G$ denote the set of all finite subsets of $G$.
 Fix an increasing  F{\o}lner sequence $\Cal F$ in $G$.
 We endow $\goth F_G$ with the discrete topology.
 Let
 $$
 \multline
 \goth R_1:=\{(C_n,F_{n-1})_{n=1}^\infty\in (\goth F_G\times\goth F_G)^{\Bbb N}\mid 
 (F_n)_{n=0}^\infty
 \text{  is a subsequence of $\Cal F$
 and}\\
\text{  \thetag{1-1}, \thetag{1-4} and \thetag{1-7}
hold}\}.
\endmultline
 $$
Denote by $\tau$ the infinite  product topology on $(\goth F_G\times\goth F_G)^{\Bbb N}$.
 Then the topological space $((\goth F_G\times\goth F_G)^{\Bbb N},\tau)$ is  Polish and  0-dimentional. 
By \cite{Da3, Lemma~3.1},  $\goth R_1$ is a $G_\delta$-subset of $(\goth F_G\times\goth F_G)^{\Bbb N}$.
Hence, $(\goth R_1,\tau)$ is a Polish space.
Define  a map $\phi:\goth R_1\to[0,+\infty]$ by setting
$\phi((C_n,F_{n-1})_{n=1}^\infty):=\lim_{n\to\infty}\frac{\#F_n}{\#C_1\cdots\#C_n}$.
Let 
 $$
  \goth R_1^{\text{fin}}:=\{\Cal T\in\goth R_1\mid\phi(\Cal T)<\infty\}.
  $$
  Of course, the condition $\phi(\Cal T)<\infty$ is equivalent to \thetag{1-3} for $\Cal T$.
  Since $\Cal F$ is F{\o}lner, it follows that \thetag{1-5} is  satisfied for each $\Cal T\in  \goth R_1^{\text{fin}}$.
  Hence, a finite measure preserving $(C,F)$-action of $G$ associated with $\Cal T$ is well defined.

We note that  $\goth R_1^{\text{fin}}$ is an $F_\sigma$-subset of  $\goth R_1$ \cite{Da3, \S3}.
 Denote by $\tau^{\text{fin}}$ the weakest topology that is stronger than $\tau$ and such that
 $\phi$ is continuous in this topology.
 Then $(\goth R_1^{\text{fin}},\tau^{\text{fin}})$ is a Polish space.
 Moreover, there is a continuous mapping $\Psi: \goth R_1^{\text{fin}}\to \Cal A_G$ such that
 $\Psi(\Cal T)$ is isomorphic to the $(C,F)$-action of $G$ associated with $\Cal T$  \cite{Da3, \S3}.
 Hence,  $\Psi(\Cal T)$ is a $G$-action of rank one along a subsequence of
 $\Cal F$.
  Conversely, each $G$-action of rank one along  a subsequence of $\Cal F$ is isomorphic to $\Psi(\Cal T)$ for some 
 $\Cal T\in \goth R_1^{\text{fin}}$ according to Fact~C.

 It was shown in \cite{Da3} that if $G$ is monotileable in the sense of \cite{We} then  the pair $(
  \goth R_1^{\text{fin}},\Psi)$ is a model for $\Cal A_G$ in the sense of \cite{Fo1}, i.e.
 for every  comeager set $M\subset \Cal A_G$ and each
$A\in M$, the set $\{\Cal T\in  \goth R_1^{\text{fin}}\mid \Psi(\Cal T)\text{ is isomorphic to $A$}\}$
 is dense in  $  \goth R_1^{\text{fin}}$.

 %$\Psi$ is one-to-one?
 %If yes then $\Psi(\Cal T)$ is a Borel subset of $\Cal A_G$.

 We let
 $$
\text{{\bf Iso}}:=\{(\Cal T,\widetilde{\Cal T})\in  \goth R_1^{\text{fin}}\times  \goth R_1^{\text{fin}}\mid \Psi(\Cal T)\text{ is isomorphic to }\Psi(\widetilde{\Cal T})\}. $$

 \proclaim{Theorem 4.1} $\text{{\bf Iso}}$ is a  $G_\delta$-subset of  $(\goth R_1^{\text{fin}}\times \goth R_1^{\text{fin}},\tau^{\text{fin}}\times \tau^{\text{fin}})$.
 \endproclaim
 
 \demo{Proof}
 Fix a decreasing sequence $(\epsilon_n)_{n=1}^\infty$ of positive reals such that  $\epsilon_1>1$ and $\sum_{n=2}^\infty\epsilon_n<0.4$.
 Given $n>0$ and a finite sequence $0=k_1<l_1<\cdots<k_n<l_n$, we say that a finite 
 sequence $ (D_m,E_{m-1}, \widetilde D_m,\widetilde E_{m-1})_{m=1}^n$
 from $(\goth F_G\times \goth F_G\times\goth F_G\times\goth F_G )^n$
 %$$
% (D_m,E_{m-1}, \widetilde D_m,\widetilde E_{m-1})_{m=1}^n
% \in (\goth F_G\times \goth F_G\times\goth F_G\times\goth F_G )^n
% $$
is {\it $(k_1,l_1,\dots,k_n,l_n)$-good} if there exist subsets
 $J_m\subset F_m$ and $\widetilde J_m\subset\widetilde F_m$  such that the following conditions are satisfied for each $m=1,\dots,n-1$:
 
 \roster
 \item"---"
 $E_m\widetilde J_m\subset\widetilde E_m$,  
 \item"---"
 $E_m^{-1}E_m\cap \widetilde J_m\widetilde J_m^{-1}=\{1\}$ 
% \item"---"
 % $\#\widetilde E_m-\# E_m\# \widetilde J_m<2\epsilon_m\#\widetilde E_m$,
 \item"---"
 $\#((\widetilde J_m J_{m+1})\triangle D_{m})<2\epsilon_m\# D_{m}$,
  \item"---"
 $\widetilde E_mJ_{m+1}\subset E_{m+1}$,   
 \item"---"
 $\widetilde E_m^{-1}\widetilde E_m\cap  J_{m+1}J_{m+1}^{-1}=\{1\}$,
% \item"---"
 %  $\# E_{m+1}-\# \widetilde E_{m}\#  J_{m+1}<2\epsilon_m\# E_{m+1}$ and
 \item"---"
 $\#(( J_{m+1}\widetilde J_{m+1})\triangle \widetilde D_{m})<2\epsilon_m\# \widetilde D_{m}$.
 \endroster
 Denote by $\Lambda(k_1,l_1,\dots,k_n,l_n)$ the subset of all 
 $(k_1,l_1,\dots,k_n,l_n)$-good sequences in $ (\goth F_G\times \goth F_G\times\goth F_G\times\goth F_G )^n$.
 Then the subset
 $$
 \multline
 V(k_1,l_1,\dots,k_n,l_n):=\{((C_n,F_{n-1})_{n=1}^\infty,( \widetilde C_n, \widetilde F_{n-1})_{n=1}^\infty)\in 
 \goth R_1^{\text{fin}}\times \goth R_1^{\text{fin}}\mid\\
 (C_{k_{m}+1,k_{m+1}}, F_{k_m},\widetilde C_{l_{m}+1,l_{m+1}}, \widetilde F_{l_m} )_{m=1}^n\text{ is good}
 \}
 \endmultline
 $$
 is clopen in $\goth R_1^{\text{fin}}$.
Moreover,  $V(k_1,l_1,\dots,k_n,l_n)$ is  $\tau$-clopen.
  Let
  $$
\Cal N:= \bigcup_{n=1}^\infty\, \bigcup_{0=k_1<l_1<\cdots<k_n<l_n}\bigcap_{l_n <k_{n+1}<l_{n+1}}V(k_1,l_1,\dots,k_n,l_n)\cap
 V(k_1,l_1,\dots,k_{n+1},l_{n+1})^c.
 $$ 
  Then $\Cal N$ is an $F_\sigma$-subset of $\goth R_1^{\text{fin}}\times \goth R_1^{\text{fin}}$.
  Hence,  the complement of $\Cal N$ is a $G_\delta$ in $\goth R_1^{\text{fin}}\times \goth R_1^{\text{fin}}$.
It follows from Theorem~2.1 that
$\Cal N^c=\text{{\bf Iso}}$.
\qed
 \enddemo
 
 \subhead 4.2. Infinite measure preserving rank-one actions
\endsubhead
In this subsection our exposition will be very sketchy. 
Let  $G$ be an arbitrary  discrete countable infinite group.
 Let $(X,\mu):=([0,+\infty), \text{Leb})$.
 Denote by Aut$(X,\mu)$ the group of $\mu$-preserving transformations of $X$.
Endow Aut$(X,\mu)$ with the weak topology, i.e. the weakest topology
in which the mappings Aut$(X,\mu)\ni R\mapsto \mu(TA\cap B)$ is continuous for all measurable subsets $A,B\subset X$ of finite measure.
Then Aut$(X,\mu)$ is a Polish group.
As in \S4.1, we denote by $\Cal A_G$ the set of all $\mu$-preserving $G$-actions on $X$.
Let
 $$
  \goth R_1^{\infty}:=\{\Cal T\in\goth R_1\mid\phi(\Cal T)=\infty\}=\goth R_1\setminus\goth R_1^{\text{fin}}.
  $$
  Of course, $\Cal T\in \goth R_1^{\infty}$ if and only if   the  $(C,F)$-action associated with $\Cal T$ is well defined and {\it infinite} measure preserving. 
On the other hand, for each infinite measure preserving rank-one $G$-action $T$, there exists $\Cal T\in\goth R_1^\infty$ such $T$ is isomorphic to the $(C,F)$-action associated with $\Cal T$.

In contrast with $\goth R_1^{\text{fin}}$, the set $ \goth R_1^{\infty}$ is a $G_\delta$-subset
of $(\goth R_1,\tau)$.
Hence, $\Cal T\in \goth R_1^{\infty}$ is a Polish space when endowed with the infinite product topology $\tau$.
Modifying slightly the construction of $\Psi$ from \cite{Da3, \S3}, one can define a continuous mapping $\Psi_\infty: \goth R_1^{\infty}\to \Cal A_G$ such that
 $\Psi_\infty(\Cal T)$ is isomorphic to the $(C,F)$-action of $G$ associated with $\Cal T$.
 Hence,  $\Psi_\infty(\Cal T)$ is a $\mu$-preserving $G$-action of rank one along a subsequence of
 $\Cal F$.
Following the argument of Theorem~4.1 almost literally, one can prove the following analogous result.

\proclaim{Theorem 4.2}
The set
$$
\text{{\bf Iso}}_\infty:=\{(\Cal T,\widetilde{\Cal T})\in  \goth R_1^{\infty}\times  \goth R_1^{\infty}\mid \Psi_\infty(\Cal T)\text{ is isomorphic to }\Psi_\infty(\widetilde{\Cal T})\}
$$
 is a  $G_\delta$-subset of   $(\goth R_1^{\infty}\times \goth R_1^{\infty},\tau\times\tau)$.
\endproclaim

 We leave details to the reader.

\enddocument

 \head Appendix
 \endhead
 
  Let $G=\Bbb Z$.
  Let $\Cal T=(C_n,F_{n-1})_{n=1}^\infty$.
 Assume that $T$ is f.m.p. and $F_n$ is $[0,1,\dots,h_n-1)$.
 Let $\epsilon_n>0$ and $\sum_{n=1}^\infty\epsilon_n<\infty$.

 \comment
 
 Suppose that there is a sequence $(\beta_n)_{n=1}^\infty$, $\beta_n>0$
 such that
 $C_{2n}+\beta_n\approx_{\epsilon_n}C_{2n}$ for each $n$.
 Then, of course, $\# C_n\to\infty$.

 We let $\alpha_n:=\beta_1+\dots+\beta_n$.
 
 Then $F_{2n}+\alpha_n\approx F_{2n}$ as $n\to
 \infty$.
 
 Define a transformation $\theta$ of $X$  by setting for each $n$:
 $$
 \theta(f_{2n}, c_{2n+1},c_{2n+2},\dots)=(\alpha_n+f_{2n}, c_{2n+1},c_{2n+2}+\beta_{n+1},c_{2n+3},c_{2n+4}+\beta_{n+2}, \dots),
 $$
 whenever this is well defined, i.e. $\alpha_n+f_{2n}\in F_{2n}$,
 $c_{2m}+\beta_{m}\in C_{2m}$ for each $m>n$. 
 It is easy to deduce from the Borel-Cantelli lemma that this transformation is well defined a.e. on $(X,\mu)$ and $\theta\in C(T)$.
 As far as I remember, I introduced $(C,F)$-transformations of such a structure in 
 
 ``{\sl Weakly mixing rank-one transformations conjugate to their squares}'',
Studia Math. 187 (2008), no. 1, 75--93.

\endcomment

 \comment
We now describe  ``the canonical structure'' of $\theta$ according to Theorem~A.
Let
$$
A_n:=C_{2n}+\alpha_n,\quad B_n:=C_{2n+1}-\alpha_n.
$$
Then we have:
$$
\align
A_n+B_n&=C_{2n}+C_{2n+1},\\
B_n+A_{n+1} &=C_{2n+1}-\alpha_n+C_{2n+2}+\alpha_{n+1}\\
&=C_{2n+1}+(C_{2n+2}+\beta_{n+1})\\
&\approx_{\epsilon_{n+1}} C_{2n+1}+C_{2n+2},\\
F_{2n-1}+A_n&\approx F_{2n},\\
F_{2n}+B_n&\approx F_{2n+1}.
\endalign
$$
Also,
let $\boldsymbol A:=(A_n)_n$ and 
$\boldsymbol B:=(B_n)_n$
and  $\psi$ be the corresponding isomorphism from $C(T)$.
Then
$$
\multline
\psi(f_{2n-1},c_{2n}, c_{2n+1},\dots)=\psi(f_{2n-1}, c_{2n}+c_{2n+1},
c_{2n+2}+c_{2n+3},\dots)\\
=\psi(f_{2n-1}, (c_{2n}+\alpha_n)+(c_{2n+1}-\alpha_n),
(c_{2n+2}+\alpha_{n+1})+(c_{2n+3}-\alpha_{n+1}),\dots)\\
=
(f_{2n-1}+(c_{2n}+\alpha_n),   (c_{2n+1}-\alpha_n)+
(c_{2n+2}+\alpha_{n+1}), (c_{2n+3}-\alpha_{n+1})+
(c_{2n+4}+\alpha_{n+2}),\dots)\\
=
(f_{2n}+\alpha_n,   c_{2n+1},
c_{2n+2}+\beta_{n+1},
 c_{2n+3},
c_{2n+4}+\beta_{n+2},\dots)
\endmultline
$$
It follows that $\theta=\psi$.
Of course, $\theta=\lim_{n\to\infty} T_{\alpha_n}$ in the weak topology.

\endcomment

%We now prove the converse. 
%The converse is equivalent to King's weak closure theorem.
Suppose that  $\theta\in C(T)$.
Then we apply our theorem on the structure of $\theta$.
Hence, after passing to an appropriate telescoping, we may assume without loss of generality that
there exist two sequences $\boldsymbol A:=(A_n)_n$ and 
$\boldsymbol B:=(B_n)_n$ of finite subsets in $G$  such that
$$
\align
F_{2n-1}+A_n &\approx F_{2n}, \  (F_{2n-1}+a)\cap (F_{2n-1}+ a')=\emptyset\text{ if }a\ne a'\tag A-1\\
F_{2n}+B_n &\approx F_{2n+1}, \  (F_{2n}+b)\cap (F_{2n}+ b')=\emptyset\text{ if }b\ne b'\tag A-2\\
A_n+B_n &\approx_{\epsilon_n}C_{2n}+C_{2n+1}\tag A-3\\
B_n+A_{n+1} &\approx C_{2n+1}+C_{2n+2}\tag A-4
\endalign
$$
for each $n$.
We let $h_n:=\# F_n$.
Then $F_n:=\{0,1,\dots, h_n-1\}$.
We note that 
$$
C_{2n}+C_{2n+1}=\bigsqcup_{c\in C_{2n+1}}(C_{2n}+c).
$$
Hence, we deduce from \thetag{A-3} that  there is $c\in C_{2n+1}$
such that
$C_{2n}+c$ is $(1-\epsilon_n)$-full of  $A_n+B_n$.
However, $A_n+B_n$ is the union of non-overlapping subsets $A_n+b$, where $b$
runs $B_n$.
The diameter of $A_n+b$ is less than $h_{2n}$.
On the other hand,
$$
\text{diam\,}(C_{2n}+c)=\text{diam}\,C_{2n}<h_{2n}.
$$
Hence $C_{2n}+c$ intersects at most two adjacent subsets $A_n+b$ and $A_n+b'$.
Thus, there is a partition $C_{2n}=C_{2n}^{(l)}\sqcup C_{2n}^{(r)}$ of $C_{2n}$ such that
$$
C_{2n}^{(l)}+c\approx (A_n+b)\cap (C_{2n}+c),\ C_{2n}^{(r)}+c\approx (A_n+b')\cap (C_{2n}+c).
$$
There is $H>0$, $H<h_{2n}$ such that
$$
F_{2n-1}+C_{2n}^{(l)}\approx [0,\dots, H)\text{ and }F_{2n-1}+C_{2n}^{(r)}\approx [H,\dots, h_{2n}).
\tag A-4
$$
%$\# A_n/ \# C_{2n}\to 1$ as $n\to\infty$.
There exist $\alpha_n$ and  $\beta_n$ such that 
$$
A_n\approx(C_{2n}^{(r)}+\alpha_n)\sqcup (C_{2n}^{(l)}+\beta_n)
$$
 and $\max( C_{2n}^{(r)}+\alpha_n)\le \min (C_{2n}^{(l)}+\beta_n)$.
%Moreover, either
%$$
%\#((C_{2n}+c)\cap(A_n+b))>(1-\epsilon_n)\#(C_{2n}+c).
%$$
%or $\#((C_{2n}+c)\cap(A_n+b'))>(1-\epsilon_n)\#(C_{2n}+c)$.
%Indeed, 
It follows from  \thetag{A-1} and \thetag{A-4} that 
$$
[H+\alpha_n, h_{2n}+\alpha_n]\sqcup [\beta_n, H+\beta_n]=[0,h_{2n}).
$$
Hence $\alpha_n=-H$, $h_{2n}-H=\beta_n$.
Thus,
$$
A_n=(C_{2n}^{(r)}-H)\sqcup (C_{2n}^{(l)}-H+h_{2n})
$$
Then $A_n+B_n$ equals
$$
(C_{2n}^{(r)}-H)\sqcup (C_{2n}^{(l)}-H+h_{2n})
\sqcup (C_{2n}^{(r)}-H+b)\sqcup (C_{2n}^{(l)}-H+h_{2n}+b)\sqcup
(C_{2n}^{(r)}-H+b')\sqcup (C_{2n}^{(l)}-H+h_{2n}+b')\sqcup\dots
$$
Then \thetag{A-3} yields that $B_n\approx C_{2n+1}+H$ and $C_{2n+1}+h_{2n}\approx C_{2n+1}$.

Then $F_{2n}+B_n\approx F_{2n+1}+H\approx F_{2n+1}$ as desired.

\enddocument